\documentclass{article}
\usepackage[a4paper,margin=1.0in]{geometry}

\usepackage{lipsum}
\usepackage{xfrac}
\usepackage{amsfonts}
\usepackage{graphicx}
\usepackage{epstopdf}
\usepackage{multirow}
\usepackage{algorithmic}
\ifpdf
  \DeclareGraphicsExtensions{.eps,.pdf,.png,.jpg}
\else
  \DeclareGraphicsExtensions{.eps}
\fi
\usepackage{subfig}
\usepackage{tikz,tikz-cd}
\usetikzlibrary{cd}
\usetikzlibrary{matrix, arrows}
\usepackage{pgfplots}
\usepackage{dsfont}
\usepackage{amsopn,bm,amssymb}
\usepackage{algorithm}
\usepackage[T3,T1]{fontenc}
\DeclareSymbolFont{tipa}{T3}{cmr}{m}{n}
\DeclareMathAccent{\invbreve}{\mathalpha}{tipa}{16}
\usepackage{amsmath}
\usepackage{amsthm}
\usepackage{cleveref}
\usepackage{authblk}

\newtheorem{theorem}{Theorem}

\newtheorem{proposition}[theorem]{Proposition}
\newtheorem{remark}[theorem]{Remark}
\newtheorem{definition}[theorem]{Definition}
\newtheorem{lemma}[theorem]{Lemma}
\newtheorem{assumption}[theorem]{Assumption}

\usepackage[normalem]{ulem}

\usepackage[colorinlistoftodos,prependcaption,textsize=tiny]{todonotes}

\usetikzlibrary{math}

\crefname{definition}{definition}{definitions}
\Crefname{definition}{Definition}{Definitions}
\Crefname{problem}{Problem}{Problems}
\crefname{problem}{problem}{problems}
\Crefname{proposition}{Prop\-o\-si\-tion}{Propositions}
\crefname{proposition}{proposition}{propositions}
\crefname{remark}{remark}{remarks}
\Crefname{remark}{Remark}{Remarks}
\crefname{algorithm}{algorithm}{algorithms}
\Crefname{algorithm}{Al\-go\-ri\-thm}{Algorithms}
\crefname{assumption}{assumption}{assumptions}
\Crefname{assumption}{Assumption}{Assumptions}
\newcommand{\cA}{\mathcal{A}}

\newcommand{\cE}{\mathcal{E}}

\newcommand{\cG}{\mathcal{G}}
\newcommand{\cH}{\mathcal{H}}

\newcommand{\cO}{\mathcal{O}}
\newcommand{\cP}{\mathcal{P}}

\newcommand{\cM}{\mathcal{M}}
\newcommand{\cN}{\mathcal{N}}

\newcommand{\cW}{\mathcal{W}}
\newcommand{\cX}{\mathcal{X}}

\newcommand{\fE}{\mathfrak{E}}
\newcommand{\fF}{\mathfrak{F}}

\newcommand{\fX}{\mathfrak{X}}

\newcommand{\IE}{\mathbb{E}}

\newcommand{\IH}{\mathbb{H}}
\newcommand{\IN}{\mathbb{N}}
\newcommand{\IP}{\mathbb{P}}
\newcommand{\IR}{\mathbb{R}}

\newcommand{\bOne}{\mathds{1}}

\renewcommand{\d}{{\text{d}}}
\newcommand{\dx}{{\d x}}
\newcommand{\dy}{{\d y}}
\newcommand{\pr}[3]{\left(#2, #3\right)_{#1}}

\newcommand{\norm}[2]{{\left\Vert #2 \right\Vert}_{#1}}

\newcommand{\abs}[1]{{\left\vert #1 \right\vert}}

\newcommand{\lp}[2]{{L}^{#1}( #2 )}

\newcommand{\probspace}{\cP, \fF, \mu}
\newcommand{\probspacesum}{\cP}
\newcommand{\probspacepar}{(\probspace )}
\newcommand{\mc}[2]{\mathtt{MC}_{#1}{\left(#2\right)}}
\newcommand{\mlmc}[2]{\mathtt{MLMC}_{#1}{\left(#2\right)}}
\newcommand{\mlmk}[2]{\mathtt{MLIPS}_{#1}{\left(#2\right)}}
\newcommand{\work}[1]{\mathtt{Wk}\!\left(#1\right)}
\theoremstyle{definition}

\DeclareMathOperator*{\esssup}{ess\,sup}

\title{A Multilevel Interacting Particle System Method for the estimation of Failure Probabilities}
\author[1]{Rub\'en Aylwin}
\author[2]{Jos\'e Pinto}
\affil[1]{Institut f\"ur Numerische Mathematik, Universit\"at Ulm}
\affil[2]{Facultad de Ingenier\'ia y Ciencias, Universidad Adolfo Ib\'añez}
\date{\today}
\begin{document}

\maketitle

% REQUIRED
\begin{abstract}
  We propose a novel method that combines a multilevel decomposition with the framework of Interacting Particle Systems to compute the probability that given quantities of interest take values above or below a certain pre-specified value. Through the application of a sequential sampling scheme, the algorithm is able to achieve convergence rates faster than those achieved by the classical Multilevel Monte Carlo algorithm, and competes with the current state of the art techniques, while also being dimension independent. Furthermore, the sample points generated by the algorithm form a sequence which concentrates near the critical value of the quantity of interest. Our results are rigorously established and later verified through numerical experiments.
\end{abstract}

% % REQUIRED
% \begin{keywords}
% Uncertainty Quantification, Multilevel Approximation, Monte Carlo Sampling, Probability Estimation.
% \end{keywords}

% % REQUIRED
% \begin{MSCcodes}
% 65C05, 65C40
% \end{MSCcodes}

\section{Introduction}
\label{sec:Intro}
During the last decades, \emph{computational uncertainty quantification} and the numerical analysis of Partial
Differential Equations (PDEs) with uncertain parameters have gained increased traction and interest in
engineering sciences and applied mathematics. The latter topic includes applications in forward uncertainty quantification,
where the distribution of certain parameters is given and one wishes to compute statistical information on given Quantities
of Interest (QoIs) and backward uncertainty quantification, where the interest lies on computing information of the
distribution of the input parameters under given values or information of QoIs. Current methods for uncertainty
quantification include: Multilevel Monte Carlo (MLMC) sampling \cite{barthMultilevelMonteCarlo2011,gilesMultilevelNestedSimulation2019},
high order Quasi Monte Carlo sampling \cite{gantnerComputationalHigherOrder2016}, high order multilevel quasi Monte Carlo sampling
\cite{dickMultilevelHigherorderQuasiMonte2017}, sparse-grid quadratures \cite{zechMultilevelApproximationParametric2019} and first-order
second-moment approaches \cite{harbrechtSparseSecondMoment2008,vonpetersdorffSparseFiniteElement2006}. Moreover, the previous algorithms
have been applied in the various engineering contexts such as: computational electromagnetics
\cite{aylwinDomainUncertaintyQuantification2020,aylwinMultilevelDomainUncertainty2023,dolzShapeUncertaintyQuantification2024,silva-oelkerQuantifyingImpactRandom2018},
linear elasticity \cite{fittUncertaintyQuantificationElastic2019,harbrechtQuantifyingDomainUncertainty2024} and acoustic scattering
\cite{dolzPMultilevelMonteCarlo2025,hiptmairLargeDeformationShape2018}, showcasing the increasing interest of the community in the
efficient estimation of uncertain responses.

A particular sub-area of uncertainty quantification that has gained additional interest corresponds to the computation of
\emph{failure probabilities}, i.e., the computation of the probability that a specific system fails (or behaves in a certain undesired
way). Failure probabilities are usually represented as the probability that a given QoI rises above (or below) a certain critical value. For examples of applications
and algorithms in this line we refer to
\cite{auEstimationSmallFailure2001,dasguptaFailureProbabilityEstimation2019,delmoralGenealogicalParticleAnalysis2005,elfversonMultilevelMonteCarlo2016,gobetRareEventSimulation2015,haji-aliAdaptiveMultilevelMonte2022,jiaDensityExtrapolationApproach2021,wagnerMultilevelSequentialImportance2020} and references therein. We particularly point to the algorithm developed in
\cite{haji-aliAdaptiveMultilevelMonte2022}, where an adaptive solver, which only considers high cost and high resolution
solutions of the PDE solver at parameter points lying close to the QoI's threshold, is used in combination with a multilevel
Monte Carlo scheme in order to compute failure probabilities with increased efficiency.

In the present article, we propose a novel multilevel Monte Carlo method, based in Interacting Particle
Systems (see \cite{cerouNonasymptoticTheoremUnnormalized2011,DelMoralFeynmanKac2004,delmoralGenealogicalParticleAnalysis2005,moralMultilevelSequentialMonte2017})
to compute the probability that a QoI takes values above (or below) a given critical value.
Through the construction of appropriate \emph{Feynman-Kac} measures, we sample the contributions of multilevel estimators only
where they might provide some information of the underlying distribution, i.e., parameters where the QoI lies close to the critical
value. This enables us to focus the sampling of the multilevel contributions only where sampling is relevant and thus decrease the number
of necessary samples on high-precision, high-cost levels in comparison to a naive multilevel Monte Carlo implementation.
Moreover, the method naturally generates samples converging to the limit surface which may be of use when computing sensitivities
  of the probability functional (see \cite{papaioannouReliabilitySensitivityEstimation2018} and references therein).

\subsection{Paper Layout}
\label{ssec:struc}
The remainder of the paper is structured as follows. In \Cref{sec:cont} we introduce the problem of computing failure probabilities
associated with functionals of solutions to parametric operator equations and their corresponding
numerical approximations. \Cref{sec:MC} is devoted to a presentation of the classical Monte Carlo and
multilevel Monte Carlo sampling strategies to our context, and is based on the expositions in
\cite{barthMultilevelMonteCarlo2011,cliffeMultilevelMonteCarlo2011,
  elfversonMultilevelMonteCarlo2016,haji-aliAdaptiveMultilevelMonte2022}.
We then construct our proposed approximation in \Cref{sec:mlips} and present proofs for its
asymptotic convergence properties. This is followed by \Cref{sec:experiments}, where we show
extensive numerical experiments, first with a synthetic example for ease of computation and
then through the consideration of a parametric Poisson equation, that our proposed algorithm
achieves the expected asymptotic convergence rates, outperforming the classical Monte Carlo and multilevel Monte Carlo strategies. We also show,
that the proposed algorithm achieves, in practice, convergence rates equivalent to those of the estimator in \cite{haji-aliAdaptiveMultilevelMonte2022}
which sets, to the best of our knowledge, the current state of the art with respect to convergence rate in this particular setting. Finally,
\Cref{sec:conclusions} contains our concluding remarks and outlines possible future research avenues.

\subsection{Notation}
\label{ssec:not}
We introduce $\IN_0:=\IN\cup \{0\}$, i.e., the naturals including $0$, $\IR^{-}:=(-\infty,0)$, and
for $a<b, \in \IN_0$ we denote the set $\{a,a+1,\hdots, b \}$ as $\{a : b\}$.
We will also employ the symbols $\lesssim$ and $\gtrsim$ to avoid specifying constants which do not depend on quantities
relevant to the corresponding analysis. 
Let $\cX$ be a given Banach space, with norm denoted by $\norm{\cX}{\cdot}$, the set of all bounded linear
functionals on $\cX$, i.e. its dual, is denoted $\cX'$. When $\cX$ is a Hilbert space, we write its
inner product as $\pr{\cX}{\cdot}{\cdot}$.

Let $\cP$ denote a sample space, $\fF$ a $\sigma$-algebra on $\cP$ and $\mu$
a probability measure on $\cP$, then we call $(\cP, \fF)$ a measurable space and
$(\cP, \fF, \mu)$ a probability space. If $(\cX,\fX)$ is another measurable space, then we call $f
:\cP\to\cX$
a measurable random variable if the preimage under $f$ of every measurable set in $\fX$
belongs to $\fF$. If $\cX$ is a Banach space, which we always consider
together with the $\sigma$-algebra generated by its open sets, then for any random variable $f:\cP\to\cX$
on the probability space $\probspacepar$, its expected value is denoted as $\IE(f)\in\cX$,
and the probability of the set $\{y\in\cP\,\vert\, f(y)\in X\}$, for $X\in\fX$, is simply written as $\IP(f\in X)$.
Though the notation for expected values and probabilities makes no direct reference to the underlying measure, this will always be clear from the context. For any event
$A\in\fF$ we write its probability as $\IP(A)$ and its complement as $A^c$. Furthermore, the
 conditional expectation of the random variable $f:\cP\to\cX$ conditioned to a measurable set $A\in\fF$
is written as $\IE(f\vert A)$, and the conditional probability of a second measurable set $B\in\fF$ conditioned
to $A$ as $\IP(B\vert A)$.

Let $\{x_n\}_{n\in\IN_0}$ be a sequence of random variables on $\probspacepar$ and taking values
on the measurable spaces $\{(\cE_n, \fE_n)\}_{n\in\IN_0}$. We say $\{x_n\}_{n\in\IN_0}$ is a Markov
chain if, for any $n>0$,
\begin{align*}
  \begin{aligned}
  &\IP\left(x_n\in A_n \middle\vert\,x_{n-1}\in A_{n-1},\hdots, x_0\in A_0\right)=\IP\left(x_n\in A_n\,\middle\vert\,x_{n-1}\in A_{n-1}\right),
    \end{aligned}
\end{align*}
where $A_n\in\fE_n,\hdots,A_0\in\fE_0$ are arbitrary. In other words, the sequence
$\{x_n\}_{n\in\IN_0}$ is a Markov chain if the state of the chain at any given time $n\in\IN_0$
depends on its past history only through its most recent state.
We further introduce, for each $n\in\IN$, the Markov
transitions $M_n:\cE_{n-1}\times\fE_n\to[0,1]$ which satisfy
$$M_n(x, A_n)=\IP\left(x_n\in A_n\middle\vert\,x_{n-1}=x\right),$$
and are, for every fixed $x_{n-1}\in\cE_{n-1}$, probability measures on $(\cE_{n}, \fE_{n})$.
Furthermore, if we take a transition $M:\cP\times\fF\to[0, 1]$ (i.e., defining a transition from $\cP$ onto itself),
we say that it has $\mu$ as an invariant measure if for every measurable set $A$ it holds that
$$\IP(A)=\int\limits_{\cP}M(x, A)\mu(\d x)=\int\limits_{\cP}\int\limits_{A}M(x, \d y)\mu(\d x),$$
i.e., the distribution induced by the probability measure $\mu$ is unaffected by the transition. For further details
on Markov chains we refer to \cite[Sec.~2.2]{DelMoralFeynmanKac2004} and \cite[Chap.~2]{hernandez-lermaMarkovChainsInvariant2003}.

For any $p\geq 1$ we denote the Bochner space of $\cX$ valued $p$-integrable random variables on
$\probspacepar$ as $\lp{p}{\probspace; \cX}$, and write its norm as $\norm{\lp{p}{\probspacesum; \cX}}{\cdot}$
when the measure and the $\sigma$-algebra are clear from the context. For more details on probability theory and Bochner spaces, we direct the interested reader to
\cite{dapratoStochasticEquationsInfinite2014,hytonenAnalysisBanachSpaces2016} and references therein.

\section{General Context}
\label{sec:cont}
\subsection{Problem Setting and Parametric Model}
\label{ssec:parpde}
Let $\cP$ be a compact subset of $\IR^s$, $s\in\IN$ and $\cH, \cM$ be two separable
Hilbert spaces. We begin by considering a parametric and uniquely solvable
operator equation given in the form
\begin{align}
  \label{eq:op}
  \cN:\cP\times\cH\to\cM'.
\end{align}
We say that $v\in\cH$ solves the operator equation given in \cref{eq:op} at a specific
point $y\in\cP$ if and only if
\begin{align}
  \label{eq:opEq}
  \cN(y, v) = 0\quad\text{in}\quad \cM',
\end{align}
and we denote the associated solution map as
\begin{align}
  \label{eq:solMap}
  u:\cP\to\cH,
\end{align}
so that the following holds:
\begin{align}
\cN(y,u(y)) = 0\quad\text{in}\quad \cM'\quad\forall y\in\cP.
\end{align}
Furthermore, we assume that the solution map remains bounded on $\cP$ so that
there is a constant $C_u>0$ such that
\begin{align}
  \label{eq:ubound}
    \norm{\cH}{u(y)}\leq C_u\quad\forall y\in\cP.
  \end{align}
Given an open subset $\cO$ of $\cH$, we wish to compute, as precisely as possible, the
probability that $u(y)$ belongs in $\cO$ when $\cP$ is granted a probability measure.
Specifically, let $\probspacepar$ be a probability space, we aim at the computation of
\begin{align}
  \label{eq:probaDefi}
  \IP(u\in\cO)=\IE\left(\bOne_{\cO}\circ u\right)=\int\limits_{\cP}\bOne_{\cO}(u(y))\mu(\dy),
\end{align}
where $\bOne_{\cO}:=\cH\to\{0,1\}$ is such that
\begin{align}
  \bOne_{\cO}(v):=
  \begin{cases}
    1\quad\text{if}\quad v\in\cO,\\
    0\quad\text{if}\quad v\not\in\cO.
  \end{cases}
\end{align}
In particular, we assume a continuous \emph{Quantity of Interest} (QoI), $g:\cH\to\IR$, bounded as
\begin{align}
  \label{eq:gbound}
  \abs{g(v)}\leq C_g\norm{\cH}{v}\quad \forall v\in\cH,
\end{align}
so that
\begin{align}
\cO:=\{v\in\cH\,\vert\, g(v) < 0\}.
\end{align}
Defining the compositions
\begin{align}
  \label{eq:qoi}
  \cG:=g\circ u:\cP\to\IR,\quad\text{and}\quad\IH_\cG:=\bOne_{\IR^-}\circ\cG,
\end{align}
we may write the probability in \cref{eq:probaDefi} as
\begin{align}
  \label{eq:probaDefi2}
  \IP(\cG<0)=\IE\left(\bOne_{\IR^-}\circ\cG\right)=\IE\left(\IH_\cG\right)=\int\limits_{\cP}\IH_\cG(y)\mu(\dy).
\end{align}
We have the following result regarding the integrability of the solution map and of the
quantity of interest.
\begin{proposition} 
  It holds for the solution map $u:\cP\to\cH$ and QoI $\cG:\cP\to\IR$ that
  \begin{align*}
    u\in\lp{2}{\cP, \fF, \mu;\cH}\quad\text{and}\quad \cG\in\lp{2}{\cP, \fF, \mu;\IR}.
  \end{align*}
\end{proposition}
\begin{proof}
  We have that
  \begin{align*}
    \abs{\cG(y)} &= \abs{g(u(y))}\leq C_g\norm{\cH}{u(y)}\leq C_gC_u,
  \end{align*}
  where $C_u,\,C_g>0$ are the constants in \cref{eq:ubound,eq:gbound} respectively,
  and the result follows from the uniform boundedness of the mappings in the statement.
\end{proof}
\subsection{Numerical Approximation}
\label{ssec:numApr}
As mentioned in the introduction, we aim at approximating the
probability in \cref{eq:probaDefi2} through a sampling scheme. Therefore, we must be able
to compute numerical approximations to the solution map in \cref{eq:solMap}.
In this context, we introduce a sequence of approximations to the
solution map depending on approximation levels $\ell\in\IN_{0}$
such that the approximation improves as the level $\ell$ increases. For any
$\ell\in\IN_{0}$ let $\cH_{\ell}$ be a finite dimensional subspace of $\cH$
such that the sequence $\{\cH_{\ell}\}_{\ell\in\IN_0}$ is dense in $\cH$.
Furthermore, we introduce
\begin{align}
  \label{eq:numSol}
  u_\ell:\cP\to\cH_\ell
\end{align}
as the \emph{level}-$\ell$ discrete solution map (which may be obtained, for example,
through the FEM, BEM or other numerical methods) and define
\begin{align}
  \label{eq:numqoi}
  \cG_\ell:=g\circ u_\ell:\cP\to\IR,\quad\text{and}\quad\IH_{\cG_\ell}:=\bOne_{\IR^{-}}\circ\cG_\ell
\end{align}
as the numerical approximations of the mappings in \cref{eq:qoi}. We make the following assumption
on the approximation properties of $\cG_\ell$.
\begin{assumption}
  \label{as:numap}
  For every $\ell\in\IN_0$ and $y\in\cP$, the computational cost of computing the
  approximation $\cG_\ell(y)=g(u_\ell(y))$ is denoted by $\cW_\ell>0$. The sequence
  $\{\cW_\ell\}_{\ell\in\IN_0}$ is monotonically increasing and there exists a constant
  $\alpha\in(0,1)$ and two exponents $q,\,r>0$ such that
  \begin{align*}
    \abs{\cG(y)-\cG_\ell(y)}\leq C_{\cG}\alpha^{q\ell}\quad\text{and}\quad \cW_\ell \lesssim\alpha^{-r\ell},
  \end{align*}
  for every $\ell\in\IN_0$ and $y\in\cP$, where both $C_{\cG}>0$ and the hidden constant are independent of both $\ell$ and $y$.
\end{assumption}
\begin{remark}
It is common that only estimates for the approximation of the solution map $u:\cP\to\cH$ are readily
available, but not for the approximation of the QoI $\cG:\cP\to\IR$, in which case the rate in
\Cref{as:numap} has to be derived from the approximation rate of the solution map through additional
assumptions. Furthermore, it could be that the exact computation of $g(u_\ell(y))$ is not
possible even for the numerical approximation $u_\ell$. In this case, we would have to consider an
additional approximation $\tilde{g}_\ell:\cH_\ell\to\IR$ of $g$ and then consider 
the convergence rate of $\tilde{g}_\ell\circ u_\ell$ to $\cG$.
\end{remark}
\section{Monte Carlo Sampling}
\label{sec:MC}
We continue our presentation by analysing the performance achieved by
straightforward implementations of single level and multilevel
Monte Carlo methods.
\subsection{Monte Carlo Sampling}
\begin{definition}
  Let $N\in\IN$, $\cX$ be a separable Hilbert space and $f:\cP\to\cX$ a measurable random variable on the probability space
  $\probspacepar$. The Monte Carlo approximation with $N$ samples of $\IE(f)$ is given by
  \begin{align}
    \mc{N}{f}:=\frac{1}{N}\sum\limits_{i=1}^Nf(y_i),
  \end{align}
  where $\{y_i\}_{i=1}^{N}$ are independent and identically distributed samples on $\cP$
  according to the density $\mu$.
\end{definition}
It is well established that if $f\in\lp{2}{\probspace; \IR}$, we have the estimate (see \cite[Lemma 4.1]{barthMultilevelMonteCarlo2011})
\begin{align}
\label{eq:erMC}
  \IE\left(\vert \IE\left(f\right)-\mc{N}{f}\vert^{2}\right)^{\frac{1}{2}}
  \leq N^{-\frac{1}{2}}\norm{\lp{2}{\cP;\IR}}{f},
\end{align}
where the outer expectation is taken with respect to the distribution of the Monte Carlo samples.
Then, for every $\ell\in\IN_0$ and $N\in\IN$, we have that
\begin{align*}
  &\IE\left(\vert \IE\left(\IH_{\cG}\right)-\mc{N}{\IH_{\cG_\ell}}\vert^{2}\right)^{\frac{1}{2}}\\
  &\leq \IE\left(\vert \IE\left(\IH_{\cG}\right)-\mc{N}{\IH_{\cG}}\vert^{2}\right)^{\frac{1}{2}}
    +
    \IE\left(\vert \mc{N}{\IH_{\cG}-\IH_{\cG_\ell}}\vert^{2}\right)^{\frac{1}{2}}\\
  &\leq N^{-\frac{1}{2}}\norm{\lp{2}{\cP;\IR}}{\IH_{\cG}}+\norm{\lp{2}{\cP;\IR}}{\IH_{\cG}-\IH_{\cG_\ell}}.
\end{align*}

We therefore require an appropriate bound for $\norm{\lp{2}{\cP;\IR}}{\IH_\cG-\IH_{\cG_\ell}}$.
The following assumption (see Assumption 2.2 in \cite{elfversonMultilevelMonteCarlo2016} and Assumption 1.4 in
\cite{haji-aliAdaptiveMultilevelMonte2022})
is required to achieve this.
\begin{assumption}
  \label{as:probreg}
  There is some $\epsilon_0>0$ such that for all positive $\epsilon < \epsilon_0$ it holds that
  \begin{align}
    \label{as:eq:probdecay}
    \IP(\cG \in (-\epsilon,\epsilon)) \lesssim \epsilon,
  \end{align}
  and it holds that $C_\cG<\epsilon_0$, where $C_\cG$ is as in \Cref{as:numap}.
\end{assumption}
\begin{remark}
  The last point in \Cref{as:probreg}, requiring $C_\cG<\epsilon_0$, enforces that the coarsest numerical approximation
  in \Cref{as:numap} (namely $\cG_0$) to be precise enough to fall under the linear decay regime of the probability in
  \cref{as:eq:probdecay}. There is no loss of generality in making this assumption, as we could dispose of it
  by stating that all the convergence results that follow hold only for sufficiently large $\ell\in\IN_0$.
\end{remark}
\begin{lemma}[Lemma 3.4 in \cite{elfversonMultilevelMonteCarlo2016}]
  \label{lem:l2proberror}
  Under \Cref{as:numap,as:probreg}, it holds that 
  \begin{align*}
    \abs{\IE(\IH_{\cG})-\IE(\IH_{\cG_\ell})}\lesssim \alpha^{q\ell}\quad\text{and}\quad
      \norm{\lp{2}{\cP;\IR}}{\IH_{\cG}-\IH_{\cG_\ell}}\lesssim \alpha^{\frac{q}{2}\ell}\quad\forall\ell\in\IN_0,
  \end{align*}
  where $q>0$ and $\alpha\in(0,1)$ are as in \Cref{as:numap}.
\end{lemma}
\begin{proof}
  Let $\ell\in\IN_0$. We have that
\begin{align*}
\abs{\IE(\IH_{\cG})-\IE(\IH_{\cG_\ell})}
  &\leq\IE\left(\abs{\IH_{\cG}-\IH_{\cG_\ell}}\right)\\
  &=\IP((\cG > 0 > \cG_\ell)\,\lor\,(\cG < 0 < \cG_\ell))\\
  &\leq\IP(\cG \in (-C_\cG\alpha^{q\ell}, C_\cG\alpha^{q\ell})),
\end{align*}
where $C_\cG>0$ is as in \Cref{as:numap}. The first part of the result then follows by employing \Cref{as:probreg}.
For the $L^2$ norm, we have that
\begin{align*}
  \norm{\lp{2}{\probspacesum;\IR}}{\IH_{\cG}-\IH_{\cG_\ell}}^2
  &=\IE\left(\abs{\IH_{\cG}-\IH_{\cG_\ell}}\right),
\end{align*}
and the result follows form the previous computation.
\end{proof}
\begin{proposition}
  \label{prop:MCRate}
  Under \Cref{as:numap,as:probreg} we can choose the number of samples $\{N_\ell\}_{\ell\in\IN_0}$ so that
  the following holds:
  \begin{align*}
    \IE\left(\vert{\IE(\IH_{\cG})-\mc{N_\ell}{\IH_{\cG_\ell}}}\vert^2\right)^{\frac{1}{2}}\lesssim\alpha^{q\ell}
    \quad\text{and}\quad
    \work{\mc{N_\ell}{\IH_{\cG_\ell}}}\lesssim \alpha^{-(2q+r)\ell},
  \end{align*}
  where the outer expected value is taken with respect to the distribution of the Monte Carlo samples,
  $\work{\mc{N_\ell}{\IH_{\cG_\ell}}}$ is the total computational work required for the computation of the
  MC approximation $\mc{N_\ell}{\IH_{\cG_\ell}}$ and the hidden constants are independent of the level $\ell\in\IN_0$.
\end{proposition}
\begin{proof}
  From our previous computations and thanks to \Cref{lem:l2proberror}, it holds that
  \begin{align*}
    \IE\left(\vert{\IE(\IH_{\cG})-\mc{N_\ell}{\IH_{\cG_\ell}}}\vert^2\right)^{\frac{1}{2}}
    \lesssim \alpha^{q\ell}+N_\ell^{-\frac{1}{2}},
  \end{align*}
  where the hidden constant is independent of the level $\ell\in\IN_0$.
  We then choose $N_\ell$ to be $\Theta(\alpha^{-2q\ell})$
  (i.e., there exist positive constants $c,C>0$ such that $c\alpha^{-2q\ell} < N_\ell < C\alpha^{-2q\ell}$)
  so that
  \begin{align*}
    \IE\left(\vert{\IE(\IH_{\cG})-\mc{N_\ell}{\IH_{\cG_\ell}}}\vert^2\right)^{\frac{1}{2}}
    \lesssim \alpha^{q\ell}.
\end{align*}
  The computational work can be bounded through \Cref{as:numap} as
  \begin{align*}
    \work{\mc{N_\ell}{\IH_{\cG_\ell}}} = N_\ell\cW_\ell\lesssim \alpha^{-2q\ell}\alpha^{-r\ell} = \alpha^{-(2q+r)\ell},
  \end{align*}
  where every hidden constant is independent of $\ell\in\IN_0$.
\end{proof}
\subsection{Multilevel Monte Carlo Method}
We continue with the multilevel Mon\-te Carlo method. Our presentation is based in the
results and analysis displayed in \cite{barthMultilevelMonteCarlo2011,cliffeMultilevelMonteCarlo2011,
  elfversonMultilevelMonteCarlo2016,haji-aliAdaptiveMultilevelMonte2022}.
\begin{definition}
  \label{defi:MLMC}
  Let $f:\cP\to\IR$ a measurable random variable on the probability space
  $\probspacepar$. Let also $\{f_\ell\}_{\ell\in\IN}$ be numerical approximations to $f$
  such that each $f_\ell:\cP\to\IR$ is a measurable random variable on $\probspacepar$ and
  $\lim\limits_{\ell\to\infty}\abs{f_\ell(y)-f(y)}=0$ for every $y\in\cP$.
  The multilevel Monte Carlo approximation at level $L\in\IN_0$ with $\cN_L:=\{N_\ell\}_{\ell=0}^L$
  samples of $\IE(f)$ is given by
  \begin{align*}
    \mlmc{\cN_L}{f}:=\sum\limits_{\ell=0}^{L}\mc{N_\ell}{f_\ell-f_{\ell-1}},
  \end{align*}
  where $f_{-1}:=0$.
\end{definition}
For $f\in\lp{2}{\probspace;\IR}$ as in \Cref{defi:MLMC} we have
\begin{align*}
  \IE\left(\vert{\IE(f_L) - \mlmc{\cN_L}{f}\vert^{2}}\right)&=\sum\limits_{\ell=0}^{L}\IE\left(\vert{\IE\left(f_\ell-f_{\ell-1}\right)-\mc{N_\ell}{f_\ell-f_{\ell-1}}}\vert^2\right)\\
  &\leq\sum\limits_{\ell=0}^{L}N_\ell^{-1}\norm{\lp{2}{\cP;\IR}}{f_\ell-f_{\ell-1}}^2,
\end{align*}
where we have employed \eqref{eq:erMC} and the independence of MC approximations at each level. The outer expected value is taken with respect to the distribution of the required samples. Moreover, if $\norm{\lp{2}{\cP;\IR}}{f_\ell-f_{\ell-1}}$
decreases quickly enough, we may take fewer samples at higher levels, achieving a better convergence rate than that in \Cref{prop:MCRate}.
\begin{lemma}
  \label{lem:l2diff}
  Under \Cref{as:numap,as:probreg} and for every $\ell\in\IN$, it holds that
  \begin{align*}
\IE\left(\vert{\IH_{\cG_\ell}-\IH_{\cG_{\ell-1}}\vert^2}\right)^{\frac{1}{2}}\lesssim\alpha^{\frac{q}{2}\ell}.
  \end{align*}
\end{lemma}
\begin{proof}
  The result follows from \Cref{lem:l2proberror} and the triangle inequality.
\end{proof}

The following result, stating the convergence properties of the MLMC estimator, is adapted to our
context from \cite{cliffeMultilevelMonteCarlo2011}, and appears also in
\cite{barthMultilevelMonteCarlo2011,elfversonMultilevelMonteCarlo2016}.

\begin{proposition}[Thm.~1 in \cite{cliffeMultilevelMonteCarlo2011}]
  \label{prop:MLMCRate}
  Under \Cref{as:numap,as:probreg} and for every $L\in\IN_0$, there is a choice for the number of samples per
  level $\{\cN_L\}_{L\in\IN_0}$ such that
  \begin{gather*}
    \IE\left(\vert{\IE(\IH_{\cG})-\mlmc{\cN_L}{\IH_{\cG}}\vert^{2}}\right)^{\frac{1}{2}}\lesssim\alpha^{qL},\\
    \work{\mlmc{\cN_L}{\IH_\cG}}\lesssim\begin{cases}
      \alpha^{-2qL}\quad &\text{if}\quad q>r\\
      L^2\alpha^{-2qL}\quad &\text{if}\quad q=r\\
      \alpha^{-(q+r)L}\quad &\text{if}\quad q<r\\
    \end{cases},
  \end{gather*}
  where the outer expectation is taken with respect to the distribution of the used samples
  and $\work{\mlmc{\cN_L}{\IH_\cG}}$ is the total computational work required for the computation the MLMC approximation and the hidden constants are independent of the level $L\in\IN_0$.
\end{proposition}

\section{Multilevel Interacting Particle System}
\label{sec:mlips}
In the present section we introduce and analyse a multilevel approximation for the failure probability
in \cref{eq:probaDefi2} based on interacting particle systems and conditional sampling. We begin in \cref{sec:MLMLIPS} by introducing the multilevel decomposition that will be used for the computation of the probabilities. Then in \Cref{ssec:IPS} we give introduction to the topic of interaction particle systems focused on the technical tools that will be used for the approximation of the multilevel decomposition.  Finally in \cref{ssec:mlmcProb} we introduce our novel approximation and the corresponding convergence results are given in \Cref{sec:MLIPSerrors}.

\subsection{Conditional Expectation of Multilevel Contributions}
\label{sec:MLMLIPS}
As we have already seen, the probability $\IP(\cG_L < 0)$ can be written, for any $L\in\IN_0$, as follows
\begin{align*}
  \IP(\cG_L < 0) = \IE(\IH_{\cG_L})=\sum\limits_{\ell=0}^{L}\IE(\IH_{\cG_\ell}-\IH_{\cG_{\ell-1}}).
\end{align*}
Recall $C_{\cG}>0$, the constant from \Cref{as:numap} and take $y\in\cP$ such that
$$\abs{\cG_{\ell}(y)}>C_{\cG}\alpha^{q\ell}(\alpha^{q}+1)$$
holds for some $\ell\in\IN_0$. Then, it follows that
\begin{align}\label{eq:integrand0}\IH_{\cG_{\ell+1}}(y)-\IH_{\cG_{\ell}}(y) = 0.\end{align}
Therefore, if we take a sequence of sets $\{\cA_\ell\}_{\ell\in\IN}$ such that
$$\{y\in\cP\,\vert\,\abs{\cG_\ell(y)}\leq C_\cG\alpha^{q\ell}(\alpha^q+1)\}\subseteq \cA_\ell$$
we will have that
\begin{align}
  \label{eq:expcond}
  \IE(\IH_{\cG_\ell}-\IH_{\cG_{\ell-1}})=\IE(\IH_{\cG_\ell}-\IH_{\cG_{\ell-1}}\vert \cA_{\ell-1})\IP(\cA_{\ell-1}).
\end{align}
Furthermore, if it holds $\cA_{\ell+1}\subset\cA_{\ell}$ for all $\ell\in\IN$, then
\begin{align}
  \label{eq:subsetsim}
  \IP(\cA_\ell):=\prod_{p=0}^\ell\IP(\cA_p\vert\cA_{p-1}),
\end{align}
with $\cA_{-1}:=\cP$. This implies that we need not sample the ML contributions on the whole parameter space $\cP$,
but rather on smaller sets $\{\cA_\ell\}_{\ell\in\IN_0}$ satisfying the conditions displayed above. The following
proposition presents a construction of the sets $\{\cA_\ell\}_{\ell\in\IN_0}$ for which these conditions hold.
\begin{proposition}
  \label{prop:setdefi}
  Let \Cref{as:numap} hold, $\{\beta_\ell\}_{\ell=0}^{L-1}\subset\IR^+$ be a finite and non-increasing sequence, and
  define the sets $\{\cA_\ell\}_{\ell=0}^{L-1}$ as follows:
  \begin{align*}
    \cA_\ell:=\left\lbrace y\in\cP\,\Bigg\vert\,\abs{\cG_\ell(y)}\leq C_\cG(\alpha^q+1)\sum\limits_{j=\ell}^{L-1}\alpha^{q j}+\beta_\ell\right\rbrace.
  \end{align*}
  Then, it holds that
  \begin{align*}
    \left\{y\in\cP\,\vert\,\abs{\cG_\ell(y)}\leq C_\cG\alpha^{q\ell}(\alpha^q+1)\right\}\subseteq\cA_\ell\quad\text{and}\quad\cA^{\ell+1}\subseteq\cA^{\ell},
  \end{align*}
  for all $\ell\in\{0:L-1\}$.
\end{proposition}
\begin{proof}
  The first property follows from the fact that $\alpha^{q\ell}\leq\sum_{l=\ell}^{L-1}\alpha^{qj}$ for all $\ell\in\{0:L-1\}$.
  Now pick any $\ell\in\{0:L-2\}$ and consider $y\in\cA_\ell^c$ such that $$\abs{\cG_{\ell}(y)}>C_\cG(\alpha^q+1)\sum_{j=\ell}^{L-1}\alpha^{qj}+\beta_\ell.$$
  Moreover, it follows directly from Assumption \ref{as:numap} that
  \begin{align*}
    \abs{\cG_{\ell}(y)}
    &\leq\abs{\cG_{\ell+1}(y)-\cG_{\ell}(y)}+\abs{\cG_{\ell+1}(y)}\\
    &\leq\abs{\cG_{\ell+1}(y)-\cG(y)}+\abs{\cG(y)-\cG_{\ell}(y)}+\abs{\cG_{\ell+1}(y)}\\
    &\leq C_\cG\alpha^{q\ell}(\alpha^{q}+1)+\abs{\cG_{\ell+1}(y)},
  \end{align*}
  which, for $y \in \mathcal{A}_\ell^c$, implies that
  \begin{align*}
    \abs{\cG_{\ell+1}(y)} &> \beta_{\ell} + C_\cG(\alpha^q+1)\sum\limits_{j=\ell+1}^{L-1}\alpha^{qj}\\
                          &\geq \beta_{\ell+1} + C_\cG(\alpha^q+1)\sum\limits_{j=\ell+1}^{L-1}\alpha^{qj}\implies y\in\cA_{\ell+1}^c,
  \end{align*}
  and we have shown that $\cA_{\ell}^c\subseteq\cA_{\ell+1}^c$, which yields the required property after taking
  the complement.
\end{proof}

By employing \Cref{prop:setdefi}, we may estimate the expectations in \cref{eq:expcond}
by sampling only on $\cA_\ell$ and disregard $\cA_\ell^c$, on which we know that the multilevel contributions
are null. This requires, however, that we generate samples from the distribution $\mu$ conditioned to $\mathcal{A}_\ell$
(\emph{cf.}~\cite[Sect.~2.1]{cerouSequentialMonteCarlo2012}), whose measure is given by
\begin{align}
  \label{eq:condMeas}
  \mu_\ell(\d y):=\frac{1}{\IP(\cA_\ell)}\bOne_{\cA_\ell}(y)\mu(\d y)\quad\forall\ell\in\IN_0.
\end{align}
Notice that by sampling from $\mu_\ell$, we are necessarily obtaining points which approximate the limit surface $\{y\in\cP\,:\,\cG(y)=0\}$. The following proposition gives a more precise notion of this observation.

\begin{proposition}
  \label{prop:accum}
  Let \Cref{as:probreg} and the conditions of \Cref{prop:setdefi} hold, with
  $$\beta_\ell:= C_\cG(\alpha^q+1)\sum\limits_{j=L}^{\infty}\alpha^{qj},$$
  so that the sets $\cA_\ell$ are chosen independently of $L\in\IN$.
  Then, for every $\ell\in\IN_0$, it holds that
  \begin{align*}
    \IP(\cA_\ell)\lesssim\alpha^{q\ell}.
  \end{align*}
  
  Moreover, if we assume $\cG:\cP\to\IR$ to be continuous, every sequence $\{y_\ell\}_{\ell\in\IN_0}\subset\cP$
  such that $y_\ell\in\cA_\ell$ for all $\ell\in\IN_0$, has a convergent subsequence whose limit point $y$
  satisfies $\cG(y)=0$.
\end{proposition}

\begin{proof}
  Let $\ell\in\IN_0$ and $y\in\cA_\ell$. Then by \Cref{as:numap} it holds that
  \begin{align*}
    \abs{\cG(y)}&\leq\abs{\cG(y)-\cG_\ell(y)}+\abs{\cG_\ell(y)}\\
                &\leq C_\cG\alpha^{q\ell} + C_\cG(\alpha^q+1)\sum\limits_{j=\ell}^{\infty}\alpha^{q j}\\
                &=C_{\cG}(\alpha^q+1)\alpha^{q\ell}\left(1+\sum\limits_{j=0}^\infty\alpha^{qj}\right).
  \end{align*}
  Then, the first claim follows from the previous computation combined with \Cref{as:probreg}.
  Let $\{y_\ell\}_{\ell\in\IN}\subset\cP$ be such that $y_\ell\in\cA_\ell$ for each $\ell\in\IN$. Since $\cP$ is compact, we have that
  there is a convergent subsequence, which we still denote $\{y_{\ell_k}\}_{k\in\IN}$, converging to a limit $y\in\cP$. Moreover, the
  continuity of $\cG$ implies that
  $$\lim\limits_{\ell\to\infty}\cG(y)-\cG(y_{\ell_k})= 0.$$
  Hence, we have
  \begin{align*}
    \abs{\cG(y)} &= \lim\limits_{k\to\infty}\abs{\cG(y_{\ell_k})}\\
                 & \leq \lim\limits_{k\to\infty}C_{\cG}(\alpha^q+1)\alpha^{q\ell_k}\left(1+\sum\limits_{j=0}^\infty\alpha^{qj}\right)=0.
  \end{align*}
\end{proof}

\subsection{Interacting Particle System}
\label{ssec:IPS}
We now turn to interacting particle systems for the generation of samples of the conditional distributions
in \cref{eq:condMeas} and the approximation of the conditional expectations and probabilities in
\cref{eq:expcond,eq:subsetsim}.
We only give a short introduction to the topic here and refer the interested reader to the references
\cite{cerouSequentialMonteCarlo2012,cerouNonasymptoticTheoremUnnormalized2011,DelMoralFeynmanKac2004,delmoralGenealogicalParticleAnalysis2005}
for further details.
 
Throughout the current section (and in accordance with the notation in
\cite{cerouSequentialMonteCarlo2012,DelMoralFeynmanKac2004}) $f$ denotes
an arbitrary measurable and bounded function on the probability space $(\probspace)$, while for each $\ell\in\IN_0$ the sets $\cA_\ell$ are assumed to be as in \Cref{prop:accum},
$G_\ell$ corresponds to a \emph{potential function}, which is also measurable on $(\probspace)$ and, without loss of generality,
assumed to be such that $0\leq G_\ell(y)\leq 1$ for all $y\in\cP$ (we will later  take $G_\ell(y)=\bOne_{\cA_{\ell}}(y)$). Moreover, let $\{x_\ell\}_{\ell\in\IN_0}$ denote a non-homogeneous
Markov Chain with transition kernels $\{{M}_{n}\}_{n\in\IN}$ from $\cP$ onto $\cP$ and assume that the transitions have $\mu$ as an invariant measure.

\begin{definition}[Feynman-Kac prediction measures]
  \label{defi:FKP}
  The \emph{Feynman-Kac} prediction
  measures associated with the tran\-si\-tions $\{{M}_{\ell}\}_{\ell\in\IN}$ and potentials
  $\{{G}_{\ell}\}_{\ell\in\IN_0}$ are denoted $\eta_\ell$ and are defined through their action on a measurable
  and bounded function $f$ as (\emph{cf.}~\cite[Def.~2.3.2]{DelMoralFeynmanKac2004})
  \begin{align*}
    \eta_\ell(f):=\frac{\IE\left(f(x_\ell)\prod\limits_{p=0}^{\ell-1}G_p(x_p)\right)}{\IE\left(\prod\limits_{p=0}^{\ell-1}G_p(x_p)\right)},
  \end{align*}
  for each $\ell\in\IN_0$. We also introduce the unnormalized Feynman-Kac prediction measure $\gamma_\ell$ through its action on $f$ as
  $$
  \gamma_\ell(f):=\IE\left(f(x_\ell)\prod\limits_{p=0}^{\ell-1}G_p(x_p)\right).
  $$
\end{definition}
Our interest in the Feynman-Kac prediction measures stems from the fact that by choosing $G_\ell(y):=\bOne_{\cA_{\ell}}(y)$ and Markov transitions $M_\ell$ that have $\mu$ as an invariant measure, one has
(\emph{cf.}~\cite[Prop.~2]{cerouSequentialMonteCarlo2012})
$$\eta_{\ell}(f)=\IE(f\vert\cA_{\ell-1})$$ 
and
$$\gamma_\ell(f)=\eta_{\ell}(f)\IP(\cA_{\ell-1})=\IE(f\vert\cA_{\ell-1})\IP(\cA_{\ell-1})$$
for each $\ell\in\IN_0$. Moreover, following \cite{delmoralGenealogicalParticleAnalysis2005},
formulas for the computation of $\IE(f\vert\cA_{\ell-1})$ and $\IE(f\vert\cA_{\ell-1})\IP(\cA_{\ell-1})$ may be recovered even for a general choice of the potentials.

For given potentials $\{{G}_{\ell}\}_{\ell\in\IN_0}$ and transitions $\{{M}_{\ell}\}_{\ell\in\IN}$, the Feynman-Kac prediction measures
$\{\eta_\ell\}_{\ell\in\IN_0}$ satisfy the following recursive equation
(\emph{cf.}~\cite[Sec.~2.5.2]{DelMoralFeynmanKac2004})
\begin{align*}
  \eta_{\ell+1}(\cdot)=\int\limits_\cP K_{\ell, \eta_{\ell}}(y, \cdot)\eta_\ell (\d y)\quad\forall\ell\in\IN_0,
\end{align*}
where, for each $\ell\in\IN_0$ we have
\begin{gather}
  \label{eq:interactionS}
  \begin{gathered}
    K_{\ell, \eta}(x,\cdot):=\int\limits_\cP M_{\ell+1}(y, \cdot)S_{\ell,\eta}(x, \dy),\\
    S_{\ell,\eta}(x,\dy) := \frac{G_\ell(y)}{\eta(G_{\ell})}\eta(\dy).
  \end{gathered}
\end{gather}
% \begin{gather}
%   K_{\ell, \eta}(x,\cdot):=\int\limits_\cP M_{\ell+1}(\hat{x}, \cdot)S_{\ell,\eta}(x, \d\hat{x}),\nonumber\\
%   S_{\ell,\eta}(x,\d\hat{x}) :=  \alpha G_\ell(x)\delta_x(\d\hat{x})+(1-\alpha G_\ell(x))\Psi_\ell(\eta)(\d\hat{x}),\label{eq:interactionS}\\
%   \Psi_\ell(\eta)(\d\hat{x}):=\frac{G_\ell(\hat{x})}{\IE_{\eta}(G_{\ell})}\eta(\d\hat{x}),\nonumber
% \end{gather} for some chosen parameter $\alpha\in[0,1]$.
A discrete $N\in\IN$ sample approximation of the measures $\{\eta_\ell\}_{\ell\in\IN_0}$,
which we denote $\{\eta^{N}_\ell\}_{\ell\in\IN}$, may be obtained through an analogous process.
We begin by generating $N$ independent and identically distributed samples according to the measure $\mu=\eta_0$,
denoted $\{x_0^{(i)}\}_{i=1}^{N}$, and set
\begin{align*}
  \eta_0^{N}(\dy):=\frac{1}{N}\sum\limits_{i=1}^N\delta_{x_0^{(i)}}(\dy).
\end{align*}
Then, for any level $\ell\in\IN_0$ and given $\{x_{\ell}^{(i)}\}_{i=1}^{N}$ (the $N$ samples of level $\ell$)
the samples for level $\ell+1$ are obtained as follows. First, we compute the intermediate points $\{\hat{x}_{\ell+1}^{(i)}\}_{i=1}^N$, which are sampled from the empirical measure\footnote{This is the discrete analogue of the transition $S(\cdot,\cdot)$ in \cref{eq:interactionS}.}
% where $\hat{x}_{\ell+1}^{(i)}$ is either chosen as ${x}_{\ell}^{(i)}$ with probability $\alpha G_{\ell}({x}_{\ell}^{(i)})$, or is
% sampled from the empirical measure
\begin{align}
  \label{eq:sampDist}
  \sum\limits_{i=1}^N\frac{G_{\ell}(x_{\ell}^{(i)})\delta_{x_{\ell}^{(i)}}(\dy)}{\sum_{j=1}^NG_{\ell}(x_{\ell}^{(j)})}.
\end{align}
The point $x_{\ell+1}^{(i)}$ is then generated from $\hat{x}_{\ell+1}^{(i)}$ by sampling from the Markov transition $M_{\ell+1}(\hat{x}_{\ell+1}^{(i)}, \dy)$,
yielding the points $\{x_{\ell+1}^{(i)}\}_{i=1}^N$, from where the discrete $N$ sample approximation of $\eta_{\ell+1}$ is defined as:
\begin{align*}
  \eta_{\ell+1}^{N}(\dy):=\frac{1}{N}\sum\limits_{i=1}^N\delta_{x_{\ell+1}^{(i)}}(\dy).
\end{align*}
\Cref{alg:sample} displays the process described above, also allowing for the generation
of $\widetilde{N}\in\{1:N\}$ samples of level $\ell+1$ (\emph{cf.}~\cite[Table~1]{moralMultilevelSequentialMonte2017}).

\begin{algorithm}
  \begin{algorithmic}[1]
    \STATE \textbf{Input}: Points $\{x^{(k)}_{\ell}\}_{k=1}^N$, Markov transition $M_{\ell+1}$, potential $G_\ell$ and number of points to be generated $\widetilde{N}\in\{1:N\}$.
    \FOR{$k\in\{1:\widetilde{N}\}$}
    \STATE Generate $\hat{x}_{\ell+1}^{(k)}$ by sampling from the empirical measure in
    \cref{eq:sampDist} as $$\hat{x}_{\ell+1}^{(k)}\sim \sum\limits_{i=1}^N\frac{G_{\ell}(x_{\ell}^{(i)})\delta_{x_{\ell}^{(i)}}(\dy)}{\sum_{j=1}^NG_{\ell}(x_{\ell}^{(j)})}.$$
    \STATE Sample ${x}_{\ell+1}^{(k)} \sim M_{\ell+1}(\hat{x}_{\ell+1}^{(k)})$.
    \ENDFOR
    \RETURN $\{{x}_{\ell+1}^{(k)}\}_{k=1}^{\widetilde{N}}$.
  \end{algorithmic}
  \caption{Generation of the samples for the approximation of $\eta_{\ell+1}$ given the samples for the approximation of $\eta_\ell$.}
  \label{alg:sample}
\end{algorithm}
The convergence properties of these approximations have been thoroughly studied, for example, in
\cite{cerouSequentialMonteCarlo2012,cerouNonasymptoticTheoremUnnormalized2011,DelMoralFeynmanKac2004}.
Our multilevel setting, however, will require that we take different numbers of samples on each level. For each $L\in\IN_0$
and non-increasing sequence $\cN_L:=\{N_\ell\}_{\ell = 0}^L$, we denote by $\eta^{\cN_L}_L$ the approximation of $\eta_L$
which is obtained following the same process as for the measures $\eta^N_L$ above, but only generating
$N_{\ell+1}\leq N_{\ell}$ samples at level $\ell+1$ from the previous $N_\ell$ samples (specifically, we only generate
$N_{\ell+1}$ intermediate points $\{\hat{x}^{(i)}_{\ell+1}\}_{i=1}^{N_{\ell+1}}$, \emph{cf.}~\Cref{alg:sample}).
Notice that the approximation depends on the number of points in each level and not only of the points in the last
level, hence we highlight the dependence on $\mathcal{N}_L$ instead of only $N_L$. We further define the corresponding
approximation of the unnormalized Feynman-Kac prediction measure as:
\begin{align*}
  \gamma_L^{\cN_L}(f):= \eta_L^{\cN_L}(f)\prod\limits_{\ell=0}^{L-1}\eta_\ell^{\cN_\ell}(G_\ell).
\end{align*}
These operators, which employ a different number of samples per level, have been
  introduced and analysed before in, for example, \cite{beskosMultilevelSequentialMonte2017,moralMultilevelSequentialMonte2017}.

\begin{remark}\label{rmk:acc}
    In \cite{DelMoralFeynmanKac2004}, it is assumed that the potential functions $\{G_\ell\}_{\ell\in\IN_0}$ are strictly positive.
    In order to consider $G_\ell:=\bOne_{\cA_\ell}$, or any other choice of potential which may be $0$ in a set of positive measure,
    an additional (accessibility) condition is considered. Namely, for every $\ell\in\IN_0$ we take $E_\ell:=\{x\in\cP\,\vert\,G_\ell(x)>0\}$
    and assume that we have that $M_{\ell+1}(x,E_{\ell+1})>0$ for every $x\in E_\ell$ together with $\mu(E_0)>0$.
    For further details we refer the interested reader to \cite[Sec.~2.4.3]{DelMoralFeynmanKac2004}.
\end{remark}

\subsection{Application to the Estimation of $\IP(\cG_L<0)$}
\label{ssec:mlmcProb}
We continue by adjusting the presentation in \Cref{ssec:IPS} to our context. For each $\ell\in\IN_0$ we take
$G_\ell(y)=\bOne_{\cA_{\ell}}(y)$. Moreover, an appropriate Markov transition kernel is constructed by first
considering an auxiliary family of transitions $\{\widetilde{M}_\ell\}_{\ell\in\IN}$ from $\cP$ onto $\cP$ and satisfying
\begin{align*}
  \widetilde{M}_{\ell}(x, \dy)\mu(\dx) = \widetilde{M}_{\ell}(y, \dx)\mu(\dy)\quad\forall x,\, y\in\cP,\, \forall\ell\in\IN.
\end{align*}
Under this condition, we have that both $\mu$ and $\mu_\ell$ are $\widetilde{M}_{\ell}$ invariant for every $\ell\in\IN$.
Following \cite{cerouSequentialMonteCarlo2012}, we introduce the transition kernels $\{M_{\ell}\}_{\ell\in\IN}$ as
\begin{align}
  \label{eq:kerDef}
  {M}_{\ell}(x, \dy) := \begin{cases}
    \delta_x(\dy)\quad &\text{if}\quad x\in \cA_{\ell-1}^c \\
    \widetilde{M}_\ell(x,\dy)\bOne_{\cA_{\ell-1}}(y) + \widetilde{M}_\ell(x,\cA_{\ell-1}^c)\delta_x(\dy)\quad &\text{if}\quad x\in \cA_{\ell-1}
  \end{cases}\quad\forall\ell\in\IN,
\end{align}
which also have $\mu$ and $\mu_\ell$ as invariant measures (\emph{cf.}~Proposition 1 in \cite{cerouSequentialMonteCarlo2012}) and satisfy
  the accessibility condition described in \Cref{rmk:acc} whenever the transitions $\widetilde{M}_\ell(x,\cdot)$ are chosen to, for example,
  assign positive measure to every open subset of $\cP$.
Under these conditions, we have
\begin{align*}
  \IE(\IH_{\cG_\ell}-\IH_{\cG_{\ell-1}}\vert\cA_{\ell-1})=\eta_\ell(\IH_{\cG_\ell}-\IH_{\cG_{\ell-1}})\quad\text{and}\quad
  \IP(\cA_\ell\vert \cA_{\ell-1})=\eta_\ell(G_{\ell}).
\end{align*}
With these properties at hand, we are ready to introduce the Multilevel Interacting Particle System (MLIPS) approximation for the probability
$\IE(\IH_{\cG})$.
\begin{definition}
  \label{def:mlmk}
  Given $L\in\IN_0$ and $\cN_L:=\{N_\ell\}_{\ell=0}^L \subset \IN$, the \emph{Multilevel Interacting Particle System} approximation  is given by
  \begin{align*}
    \mlmk{\cN_L}{\IH_\cG}:=\sum\limits_{\ell=0}^{L}\gamma_\ell^{\cN_\ell}(\IH_{\cG_\ell}-\IH_{\cG_{\ell-1}}).
  \end{align*}
\end{definition}

The operator introduced in \Cref{def:mlmk} is computed in practice through \Cref{alg:comp}.

\begin{remark}
  The observant reader will notice that during the call to \Cref{alg:sample} in line 17 of \Cref{alg:comp}, we repeat
  (expensive) evaluations of the potentials $G_\ell$. It is clear that an efficient implementation will perform these evaluations
  only once and then save the results for later usage. Similarly, in line 9 of \Cref{alg:sample} (and given our specific choice of
  Markov transition in \cref{eq:kerDef}) we perform evaluations of the potential that should also be saved for their later usage in
  line 7 of \Cref{alg:comp}. We refer to \cite{githubCodeC} for an efficient Python implementation of \Cref{alg:comp}.
\end{remark}

\begin{algorithm}[ht!]
  \begin{algorithmic}[1]
    \STATE \textbf{Input}: level $L$, convergence rate $q>0$, constant $C_\cG>0$, positive non-increasing sequence $\{\beta_\ell\}_{\ell=0}^{L-1}$, samples per level $\cN_L:=\{N_\ell\}_{\ell=0}^L$, measure $\mu$, Markov transition kernels $\{M_\ell\}_{\ell=1}^L$.
    \STATE Generate $N_0$ i.i.d.~samples according to $\mu$ and denote them as $S_0:=\{x_0^{(i)}\}_{i=1}^{N_0}$.
    \STATE Set $\mathsf{MLIPS}=0$, $p=1$.
    \FOR{$\ell\in\{0:L\}$}
    \STATE Set $R_\ell = 0$ and $p_\ell = 0$.
    \FOR{$i\in\{1:N_\ell\}$}
    \STATE Compute $\cG_\ell(x^{(i)}_\ell)$ and $\cG_{\ell-1}(x^{(i)}_\ell)$.
    \STATE Compute the update $R_\ell = R_\ell + \IH_{\cG_\ell}(x^{(i)}_\ell)-\IH_{\cG_{\ell-1}}(x^{(i)}_\ell).$
  \IF{$\vert{\cG_\ell(x^{(i)}_\ell)}\vert<C_\cG\alpha^{q+1}\sum^{\infty}_{j=\ell}\alpha^{qj}$}
  \STATE Update $p_\ell = p_\ell + 1$.
  \ENDIF
  \ENDFOR
  \STATE Compute the update
  \begin{align*}
    % \label{eq:alg:evaluation}
    \begin{gathered}
      \mathsf{MLIPS}=\mathsf{MLIPS}+p\frac{R_\ell}{N_\ell}\quad\text{and}\quad
      p=p\frac{p_\ell}{ N_\ell}.
    \end{gathered}
  \end{align*}
  \IF{$\ell=L$}
  \RETURN $\mathsf{MLIPS}$.
  \ENDIF
  \STATE Call \Cref{alg:sample} with the following input: points $S_{\ell}$, Markov transition $M_{\ell+1}$, potential $G_\ell:=\bOne_{\cA_\ell}$
  and number of points to be generated $N_{\ell+1}$. Obtain $S_{\ell+1}:=\{x_{\ell+1}^{(i)}\}_{i=1}^{N_{\ell+1}}$.
  \ENDFOR
\end{algorithmic}
\caption{Algorithm for the computation of the MLIPS estimator.}
\label{alg:comp}
\end{algorithm}

\subsection{Computation and Error Analysis of the MLIPS Approximation}
\label{sec:MLIPSerrors}
The following section is concerned with the properties of the approximation introduced in
\Cref{def:mlmk}. We begin by presenting the approximation properties for the discrete versions of the
Feynman-Kac measures, which follow from straightforward modifications to the convergence results for
the approximate measures $\eta_\ell^N$. We require, however, some additional conditions on
the Markov transition kernels (see assumptions (A1) and (A2) in \cite{beskosMultilevelSequentialMonte2017,moralMultilevelSequentialMonte2017}).

\begin{assumption}
  \label{as:Mm}
  There is some $c\in(0,1)$ such that for any $p\in\IN$ and $x_p$, $y_p$ in $\cA_p$ it holds that
  \begin{align*}
    M_{p+1}(x_p, \cA) \geq c M_{p+1}(y_p, \cA),
  \end{align*}
  for any $x_p, y_p\in\cA_p$ and measurable subset $\cA\subseteq\cP$.
  Moreover, for any $p\in\IN$ and $x\in\cA_{p-1}$ it holds that $M_p(x,\cA_p)>0$.
\end{assumption}

% \begin{lemma}[Thm.~7.4.4 in \cite{DelMoralFeynmanKac2004}]
%   \label{lem:convSimp1}
%     Let \Cref{as:numap,as:probreg,as:Mm} hold. Then, for any measurable $f:\cP\to\IR$ and non-increasing
%   sequence $\cN_\ell:=\{N_p\}_{p=0}^\ell\subset\IN$, it holds that
%   \begin{align*}
%     \IE\left(\vert\IE_{\eta^{\cN_\ell}_\ell}(f)-\eta_\ell(f)\vert^2 \right)^\frac{1}{2}
%     \lesssim N_\ell^{-\frac{1}{2}}\sup\limits_{x,y\in\cP}\abs{f(x)-f(y)},
%   \end{align*}
%   where the expectation is taken over the distribution of the samples $\{x^{(k)}_{\ell}\}_{k=1}^{N_\ell}$ associated
%   with the discrete measure $\eta_\ell^{\cN_\ell}$ and the hidden constants are independent of $\ell\in\IN_0$ and $f$. 
% \end{lemma}
% \begin{proof}
%   The result follows form a straightforward modification of the proof of Theorem 7.4.4 in \cite{DelMoralFeynmanKac2004}. A short discussion on the outline of the proof is given in \Cref{apdx:outline1}.
% \end{proof}

\begin{lemma}[Proposition B.2 in \cite{moralMultilevelSequentialMonte2017}]
  \label{lem:convSimp}
  Let \Cref{as:numap,as:probreg,as:Mm} hold. Then, for any $\ell\in\IN_0$ and any given non-increasing
  sequence $\cN_\ell:=\{N_p\}_{p=0}^\ell\subset\IN$ and measurable and bounded function $f$
  \begin{align*}
    \IE\left(\vert\gamma^{\cN_\ell}_\ell(f)-\gamma_\ell(f)\vert^2 \right)^\frac{1}{2}
    \lesssim N_\ell^{-\frac{1}{2}}(\ell+1)^{\frac{1}{2}}\gamma_\ell(1)\norm{\lp{\infty}{\cP;\IR}}{f},
  \end{align*}
  where the expectation is taken over the distribution of the sample associated
  with the discrete measures $\{\eta_\ell^{\cN_\ell}\}_{\ell=0}^L$, and the hidden constants are independent of $\ell\in\IN_0$ and $f$. 
\end{lemma}
\begin{proof}
  The lemma is a direct consequence of Proposition B.2 in \cite{moralMultilevelSequentialMonte2017}
  (under Assumptions (A1) and (A2) in \cite{beskosMultilevelSequentialMonte2017,moralMultilevelSequentialMonte2017})
  and our choice of (decreasing) samples per level.
\end{proof}

% \begin{remark}
%   \ra{
%     In \cite{DelMoralFeynmanKac2004}, the analysis is also carried out under the assumption that for all $p$ and $x_p$, $y_p\in\cP$, it
%     holds that $G_p(x_p)\geq \varepsilon_p G_p(y_p) > 0$ for some $\varepsilon_p\in (0,1]$, which precludes potential functions which
%     are $0$ on sets of positive measure. As mentioned throughout \cite{DelMoralFeynmanKac2004}, this is remedied by considering the restriction
%     of the Markov transition kernels to the sets where the potentials $G_p$ are strictly positive together with the accessibility condition
%     $(\cA)$ introduced in \cite[Sec.~2.4.3]{DelMoralFeynmanKac2004}, stating that for all $x\in\cP$, $M_p(x,\{y\in\cP\,\vert\, G_p(y)>0\})>0$,
%     implying that the sets $\{y\in\cP\,\vert\, G_p(y)>0\}$ are accessible from any point in $\cP$.}
% \end{remark}

\begin{theorem}
  \label{thm:main}
  Let all the assumptions of \Cref{lem:convSimp} be satisfied. Then, for every $L\in\IN_0$, there is a choice
  for the number of samples per level $\{\cN_L\}_{L\in\IN_0}$, such that
  \begin{align*}
    \IE\left(\vert\IE(\IH_{\cG})-\mlmk{\cN_L}{\IH_{\cG}}\vert^2\right)^{\frac{1}{2}}\lesssim\alpha^{qL},
  \end{align*}
  where the hidden constant does not depend on $L\in\IN_0$. Moreover, the total work required to compute the MLIPS approximation is bounded as
  \begin{align*}
    \work{\mlmk{\cN_L}{\IH_\cG}}\lesssim\begin{cases}
      \alpha^{-2qL}\quad\text{if}\quad q>\frac{1}{2}r\\
      L^4\alpha^{-2qL}\quad\text{if}\quad q=\frac{1}{2}r\\
      L\alpha^{-rL}\quad\text{if}\quad q<\frac{1}{2}r\\
    \end{cases}.
  \end{align*}
  
\end{theorem}
\begin{proof}
  Take arbitrary $\ell\in\{1:L\}$ and introduce, for ease of notation, $f_\ell:=\IH_{\cG_\ell}-\IH_{\cG_{\ell-1}}$. Then, by \Cref{lem:convSimp} we have that
  \begin{align*}
    \IE\left(\vert\gamma^{\cN_\ell}_\ell(f_\ell)-\gamma_\ell(f_\ell)\vert^2 \right)^\frac{1}{2}
    &\lesssim N_\ell^{-\frac{1}{2}}(\ell+1)^{\frac{1}{2}}\gamma_\ell(1)\\
    &= N_\ell^{-\frac{1}{2}}(\ell+1)^{\frac{1}{2}}\IP(\cA_{\ell-1})
      \lesssim N_\ell^{-\frac{1}{2}}(\ell+1)^{\frac{1}{2}}\alpha^{q\ell}
    % &\norm{\lp{2}{\probspacesum; \IR}}{\eta_\ell(f_\ell)\IP(\cA_{\ell-1})-\IE_{\eta_\ell^{\cN_\ell}}(f_\ell)\prod\limits_{p=0}^{\ell-1}\IE_{\eta_p^{\cN_p}}(G_p)}\\
    % &\hspace{.5cm}\leq\norm{\lp{2}{\probspacesum; \IR}}{\IP(\cA_{\ell-1})-\prod\limits_{p=0}^{\ell-1}\IE_{\eta_p^{\cN_p}}(G_p)}
    % +\norm{\lp{2}{\probspacesum; \IR}}{\eta_\ell(f_\ell)-\IE_{\eta_\ell^{\cN_\ell}}(f_\ell)}\IP(\cA_{\ell-1})\\
    % &\hspace{.5cm}\lesssim \ell^{\frac{1}{2}}\IP(\cA_{\ell-1})N_{\ell-1}^{-\frac{1}{2}}+\IP(\cA_{\ell-1})N_\ell^{-\frac{1}{2}}\lesssim
    %   (\ell^{\frac{1}{2}}+1)\alpha^{q\ell}N_{\ell}^{-\frac{1}{2}},
  \end{align*}
  where all hidden constants are independent of the level $\ell\in\{1:L\}$ and the last inequality follows from \Cref{prop:accum}. For $\ell=0$ we have the usual MC error bound, leading us to
  \begin{align*}
    \IE\left(\vert\IE(\IH_{\cG})-\mlmk{\cN_L}{\IH_{\cG}}\vert^2\right)^{\frac{1}{2}}&\lesssim\alpha^{qL}+
    \sum\limits_{\ell=0}^L(\ell+1)^{\frac{1}{2}}\alpha^{q\ell}N_{\ell}^{-\frac{1}{2}},
  \end{align*}
  and the remainder of the proof follows by a simple modification of that of \cite[Thm.~1]{cliffeMultilevelMonteCarlo2011} which we defer to \Cref{apdx:samp}.
\end{proof}

\begin{remark}
  The convergence results exposed in \Cref{thm:main} are an evident improvement when compared to
  those in \Cref{prop:MLMCRate} and \Cref{prop:MCRate}. However, in
  \cite{elfversonMultilevelMonteCarlo2016} a different way of improving the results of the MLMC
  method is presented which, leveraging similar ideas to those employed here, allows for a reduction
  in the cost of computing the estimator and results in a similar (if not slightly better)
  convergence than that of the MLIPS method. Specifically, their convergence is equal for the case
  $q>\frac{1}{2}r$ while the cost is increased in a factor $L$ for $q\geq\frac{
    1}{2}r$. The method
  presented in \cite{elfversonMultilevelMonteCarlo2016}, however, can not ensure that any number of
  samples will be recognized as being close to the limit surface $\{y\in\cP\,\vert\,\cG(y)=0\}$.
\end{remark}

\begin{remark}
As mentioned in \cite[Sec.~4.3]{cerouSequentialMonteCarlo2012}, the mixing properties of the transition
$\widetilde{M}_\ell$ directly affect the performance of \Cref{alg:comp}. If the transitions are too large, then
the majority of the generated points will be rejected, which will result in a poor sample for the coming level in the estimation.
The same will hold true if the transition is too small and the generated points are mostly accepted but too
close to each other. We follow the suggestion in \cite[Sec.~4.3]{cerouSequentialMonteCarlo2012} and consider
modifying the transition kernel after each application on the sample set if it has an acceptance rate below or above
certain parameters. The mixing properties of the Transition $M_\ell$ can be further increased by repeated application
of the transition $\widetilde{M}_\ell$, though at an increased cost. In our experiments in \Cref{sec:experiments}
(with implementation available in \cite{githubCodeC}), the transition $M_\ell$ consists of three sequential applications of
$\widetilde{M}_\ell$.
\end{remark}

\section{Experimental results}
\label{sec:experiments}
We continue by presenting experimental results de\-mon\-stra\-ting the properties of the presented algorithm.
Our implementation was carried out in Python and is available in \cite{githubCodeC}. Throughout our experiments we
take $\cP=[-1,1]^s$ for some $s\in\IN$ (different for each experiment) and choose, for each $\ell\in\IN$,
the measures $\widetilde{M}_\ell(x,\d y)$ to follow a Normal distribution centered at
$x$ and with a diagonal covariance matrix $\Sigma_\ell$ and restricted to the interval $\cP$ by
performing a repeated reflection along the boundaries of $\cP$. Then, since the densities are bounded from below
and above, the conditions in \Cref{lem:convSimp} are clearly satisfied.
\subsection{A Simple Demonstrative Example}
To showcase the properties of \Cref{alg:comp} we consider, first, the following synthetic example.
We take $\cP:=[-1,1]^2$, $\mu$ as the normalized Lebesgue measure on $\cP$ and define, for $y:=(y_1,y_2)\in\cP$,
\begin{align*}
  \cG(y):=\frac{1}{4}(y_1^2+y_2^2)
  \quad\text{and}\quad
  \cG_\ell(y):=\cG(y)+\varepsilon 2^{-q\ell}h_\ell(y) \forall\ell\in\IN_0,
\end{align*}
where $\varepsilon>0$ and, for every $\ell\in\IN_0$ it holds that
\begin{align*}
  h_\ell\in\lp{2}{\cP}
  \quad\text{and}\quad
  \esssup\limits_{y\in\cP}\abs{h_\ell(y)}=1,
\end{align*}
so that \Cref{as:numap} is satisfied with constant $\varepsilon$ for the error of the numerical method and
$\alpha=\frac{1}{2}$. Specifically, we shall take $h_\ell(y):=\tfrac{1}{2}(\sin(\tfrac{1}{\ell+1}\pi y_1)-\cos(\tfrac{1}{\ell+1}\pi y_2))$.
Furthermore, for any $\vartheta\in(0,\frac{1}{4})$ we have that
$\IP(\cG(y)<\vartheta)=\pi\vartheta$. Additionally, we take the work at level $\ell\in\IN_0$ to be equal
to $2^{r\ell}$, and we shall perform several numerical experiments with varying values of the rates $q,r>0$. We shall also take different
values of the cutoff $\vartheta$ ($0.1$ and $0.01$), to study the performance of the algorithms for events with different probabilities,
the error constant $\epsilon$ (namely $0.1\vartheta$ and $0.05\vartheta$), to simulate different base precisions of the numerical method.
The error constant $\epsilon$ is chosen as a fraction of $\vartheta$ so as to reflect that the resolution of the chosen numerical methods 
should depend on their behavior near the cutoff $\vartheta$. The convergence diagrams that will
be presented in the following sections display the relative error as a function of the expended computational work.

\subsubsection{$q>\frac{1}{2}r$}
\label{sssec:numqgr}
\begin{figure*}[ht!]
  \subfloat[$\vartheta=0.1$ and $\epsilon=0.05\vartheta$]{%
    \includegraphics[width=.4\linewidth]{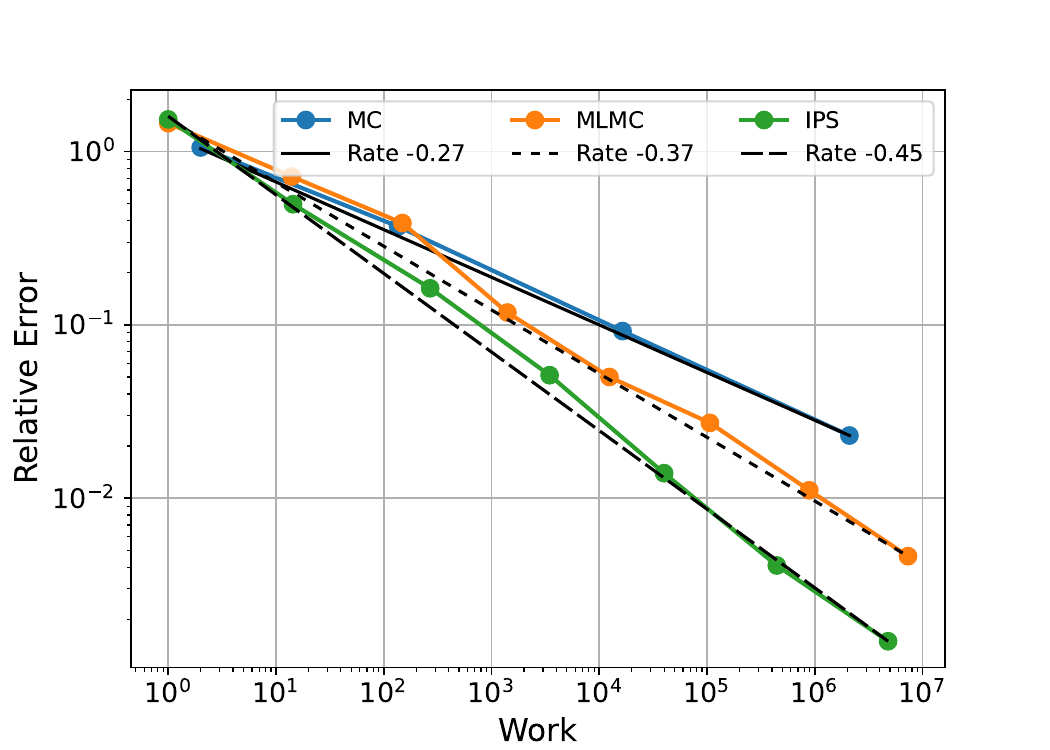}%
    \label{subfig:q2r3:a}%
  }\hfill
  \subfloat[$\vartheta=0.1$ and $\epsilon=0.01\vartheta$]{%
    \includegraphics[width=.4\linewidth]{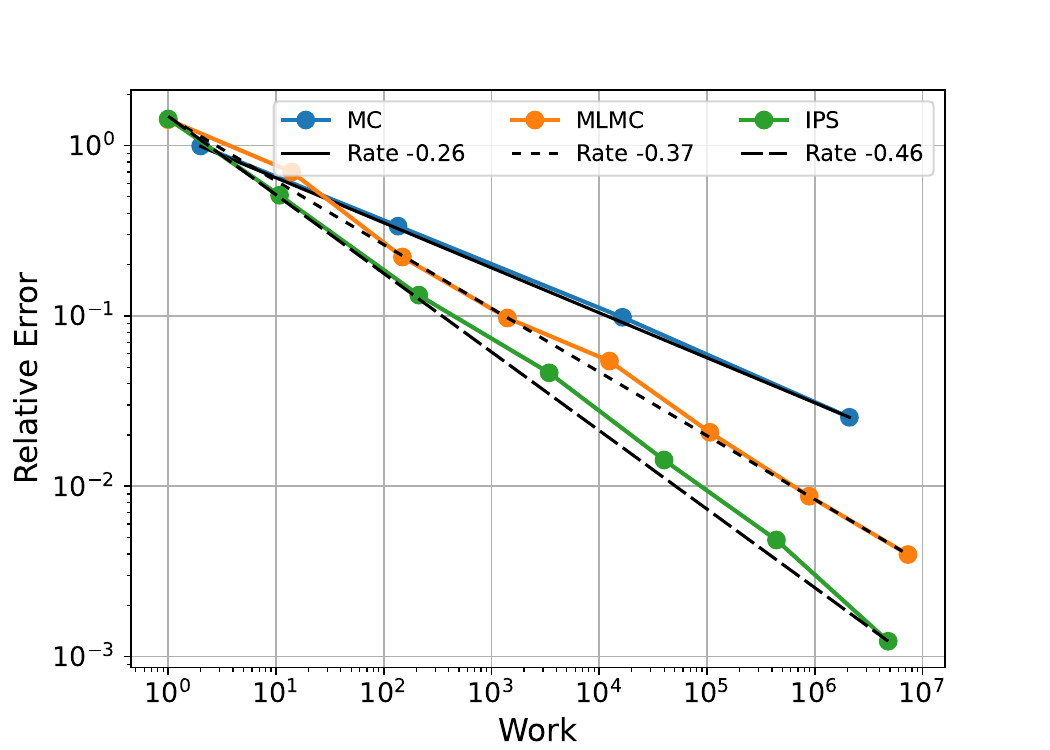}%
    \label{subfig:q2r3:b}%
  }\\
  \subfloat[$\vartheta=0.01$ and $\epsilon=0.05\vartheta$]{%
    \includegraphics[width=.4\linewidth]{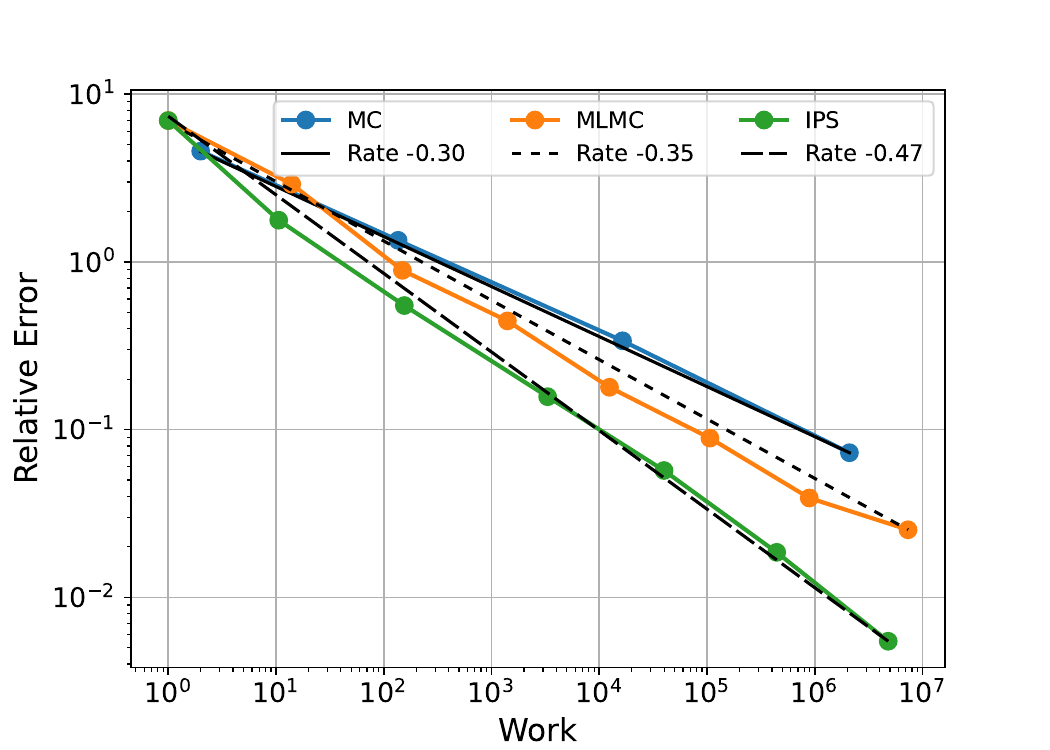}%
    \label{subfig:q2r3:c}%
  }\hfill
  \subfloat[$\vartheta=0.01$ and $\epsilon=0.01\vartheta$]{%
    \includegraphics[width=.4\linewidth]{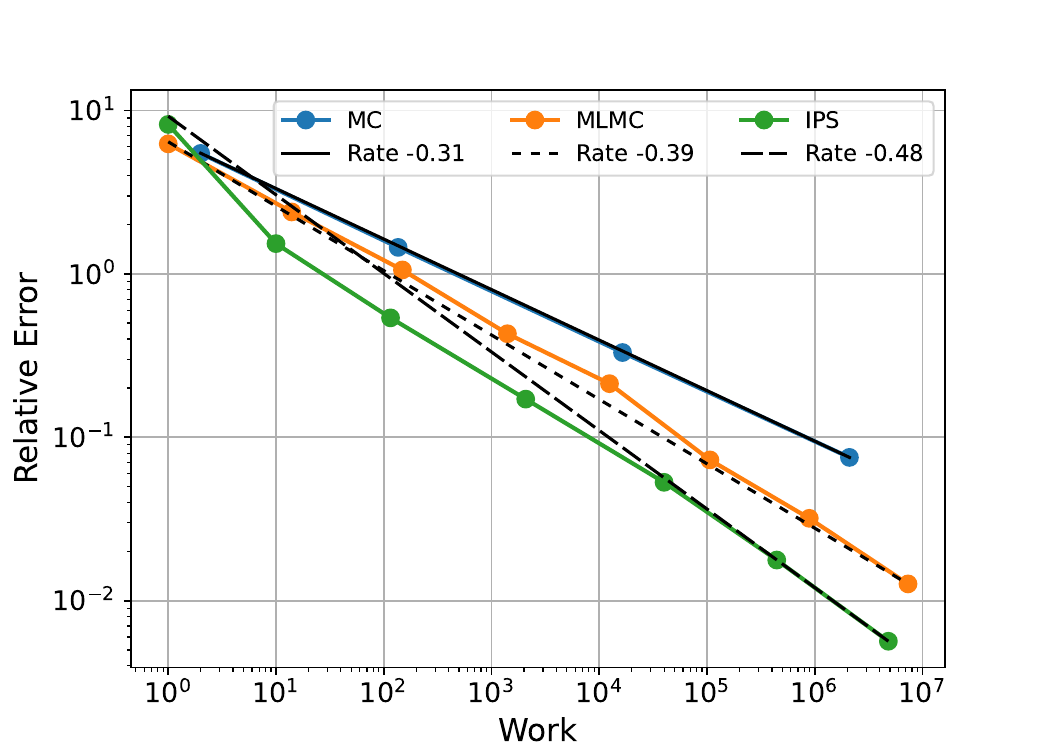}%
    \label{subfig:q2r3:d}%
  }
  \caption{Convergence for the Monte Carlo method, multilevel Monte Carlo method and multilevel interacting particle system method for $q=2$ and $r=3$ for different values of the cutoff $\vartheta$ and error constant $\epsilon$.}
  \label{fig:Simpleq2r3}
  % Run with:
  % python runExperiment.py --ips --mlmc --mc -L 7 -Y 0.1 -E 0.005 -r 3 -q 2 -g 0.5 --realizations 200 -f E.005_Y0.1_r3_q2_g0.5
  % python runExperiment.py --ips --mlmc --mc -L 7 -Y 0.1 -E 0.001 -r 3 -q 2 -g 0.5 --realizations 200 -f E.001_Y0.1_r3_q2_g0.5
  
  % python runExperiment.py --ips --mlmc --mc -L 7 -Y 0.01 -E 0.0005 -r 3 -q 2 -g 0.5 --realizations 200 -f E.0005_Y0.01_r3_q2_g0.5
  % python runExperiment.py --ips --mlmc --mc -L 7 -Y 0.01 -E 0.0001 -r 3 -q 2 -g 0.5 --realizations 200 -f E.0001_Y0.01_r3_q2_g0.5
\end{figure*}

We begin with the case $q>\frac{1}{2}r$. We also impose $q<r$ to showcase the improved convergence
achieved by the MLIPS method when compared to the MLMC method. Therefore, we take $q=2$ and $r=3$, values for which we expect convergence rates of: $\frac{q}{2q+r}=\frac{2}{7}$ for the Monte Carlo method, $\frac{q}{q+r}=\frac{2}{5}$ for the multilevel Monte Carlo method, and $\frac{q}{2q}=\frac{1}{2}$ for the multilevel interacting particle system method.
\Cref{fig:Simpleq2r3} displays the convergence rate of the method presented in the previous sections for
this case, where each point represents the average squared error of $100$ realizations. We observe the MC method achieving its expected convergence rate for all cases. The same holds true for the MLMC method, while the MLIPS method achieves convergence rates which are close to the theoretical bound of $\frac{1}{2}$. Moreover, we see an upwards shift in the relative error as the probability we are computing decreases, displaying the dependence of the error bound of all these methods on the failure probability they aim to compute (compare the results in \cref{subfig:q2r3:a,subfig:q2r3:b} with those in \cref{subfig:q2r3:c,subfig:q2r3:d}, where the probability of interest is $10$ times smaller).

\subsubsection{$q=\frac{1}{2}r$}
\begin{figure*}[ht!]
  \subfloat[$\vartheta=0.1$ and $\epsilon=0.05\vartheta$]{%
    \includegraphics[width=.4\linewidth]{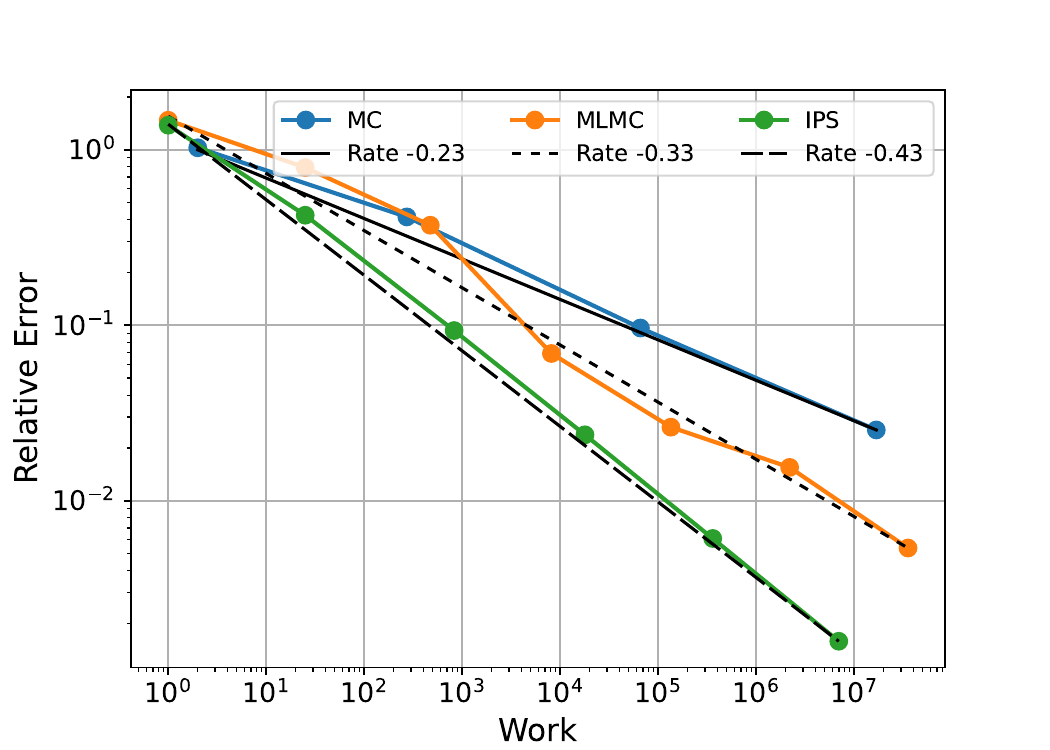}%
    \label{subfig:q2r4:a}%
  }\hfill
  \subfloat[$\vartheta=0.1$ and $\epsilon=0.01\vartheta$]{%
    \includegraphics[width=.4\linewidth]{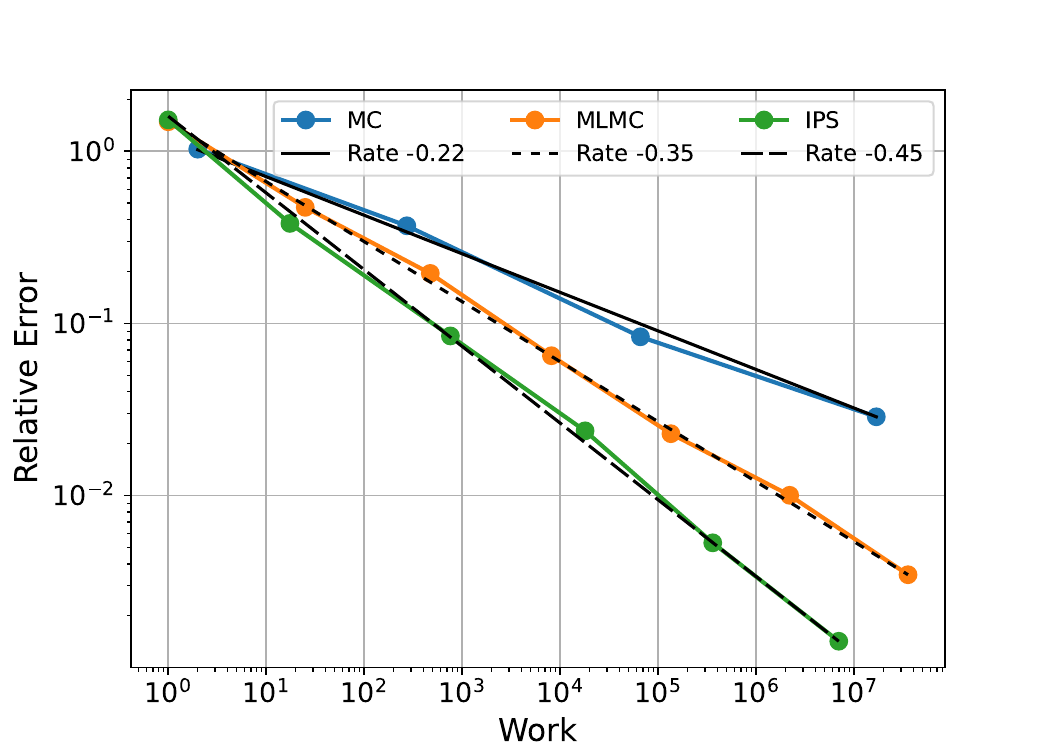}%
    \label{subfig:q2r4:b}%
  }\\
  \subfloat[$\vartheta=0.01$ and $\epsilon=0.05\vartheta$]{%
    \includegraphics[width=.4\linewidth]{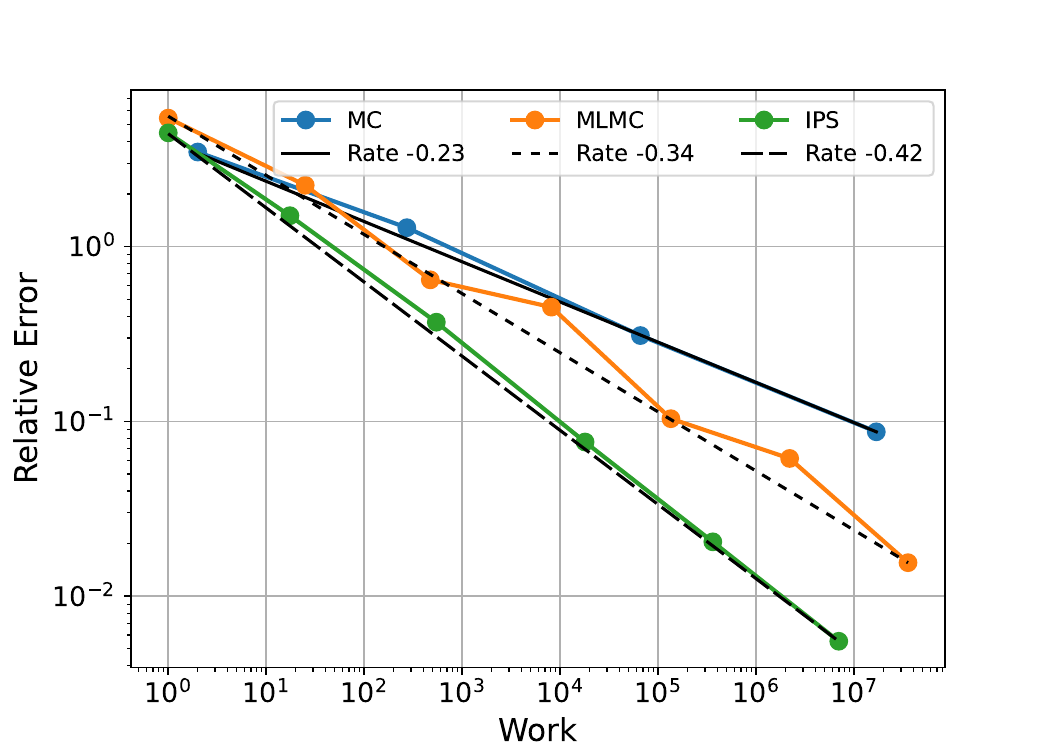}%
    \label{subfig:q2r4:c}%
  }\hfill
  \subfloat[$\vartheta=0.01$ and $\epsilon=0.01\vartheta$]{%
    \includegraphics[width=.4\linewidth]{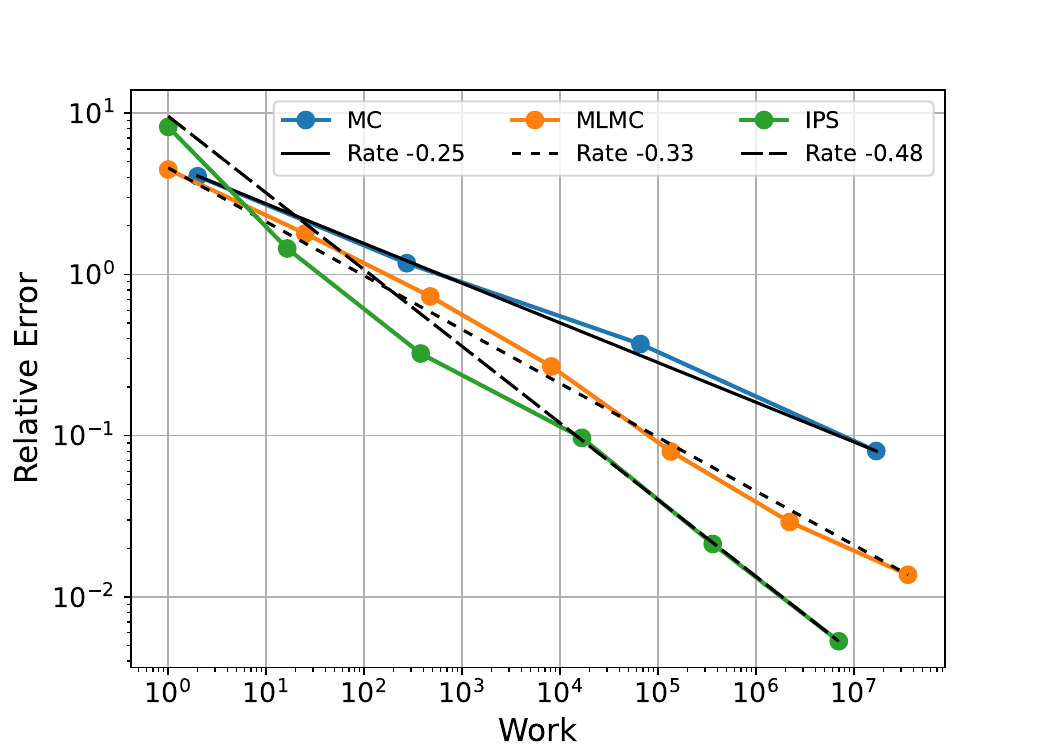}%
    \label{subfig:q2r4:d}%
  }
  \caption{Convergence for the Monte Carlo method, multilevel Monte Carlo method and multilevel interacting particle system method for $q=2$ and $r=4$ for different values of the cutoff $\vartheta$ and error constant $\epsilon$.}
  \label{fig:Simpleq2r4}
  % Run with:
  % python runExperiment.py --ips --mlmc --mc -L 8 -Y 0.1 -E 0.005 -r 4 -q 2 -g 0.5 --realizations 100 -f E.005_Y0.1_r4_q2_g0.5
  % python runExperiment.py --ips --mlmc --mc -L 8 -Y 0.1 -E 0.001 -r 4 -q 2 -g 0.5 --realizations 100 -f E.001_Y0.1_r4_q2_g0.5
  
  % python runExperiment.py --ips --mlmc --mc -L 8 -Y 0.01 -E 0.0005 -r 4 -q 2 -g 0.5 --realizations 100 -f E.0005_Y0.01_r4_q2_g0.5
  % python runExperiment.py --ips --mlmc --mc -L 8 -Y 0.01 -E 0.0001 -r 4 -q 2 -g 0.5 --realizations 100 -f E.0001_Y0.01_r4_q2_g0.5
\end{figure*}

We now take $q=2$ and $r=4$ and repeat our previous experiment, which are presented in \Cref{fig:Simpleq2r4}, for which we expect convergence rates of $\frac{q}{2q+r}=\frac{1}{4}$ for the MC method and $\frac{q}{q+r}=\frac{1}{3}$ for the MLMC. For the MLIPS method we no longer expect a rate close to $\frac{q}{2q}=\frac{1}{2}$, but a somewhat diminished convergence rate due to the presence of the term $L^4$ in the estimate shown in \Cref{thm:main} for this case. \Cref{fig:Simpleq2r4} displays the obtained convergence results for the present case, where we observe a similar behavior as for the previous case, with all methods achieving their expected theoretical convergence rate, with only a small worsening of the convergence rate of the multilevel interacting particle system method in comparison to that in the previous case ($q>\frac{1}{2}r$). Once again, we observe higher relative errors for the experiments performed for a smaller value of the cutoff $\vartheta$.

\subsubsection{$q<\frac{1}{2}r$}
\begin{figure*}[ht!]
  \subfloat[$\vartheta=0.1$ and $\epsilon=0.05\vartheta$]{%
    \includegraphics[width=.4\linewidth]{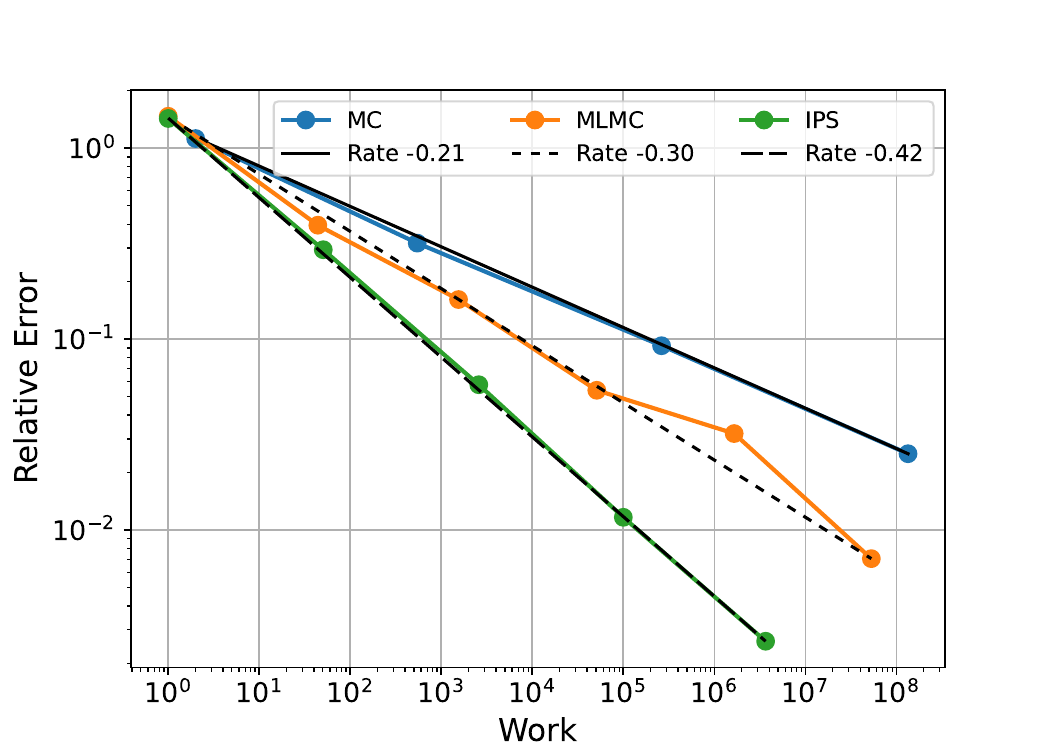}%
    \label{subfig:q2r5:a}%
  }\hfill
  \subfloat[$\vartheta=0.1$ and $\epsilon=0.01\vartheta$]{%
    \includegraphics[width=.4\linewidth]{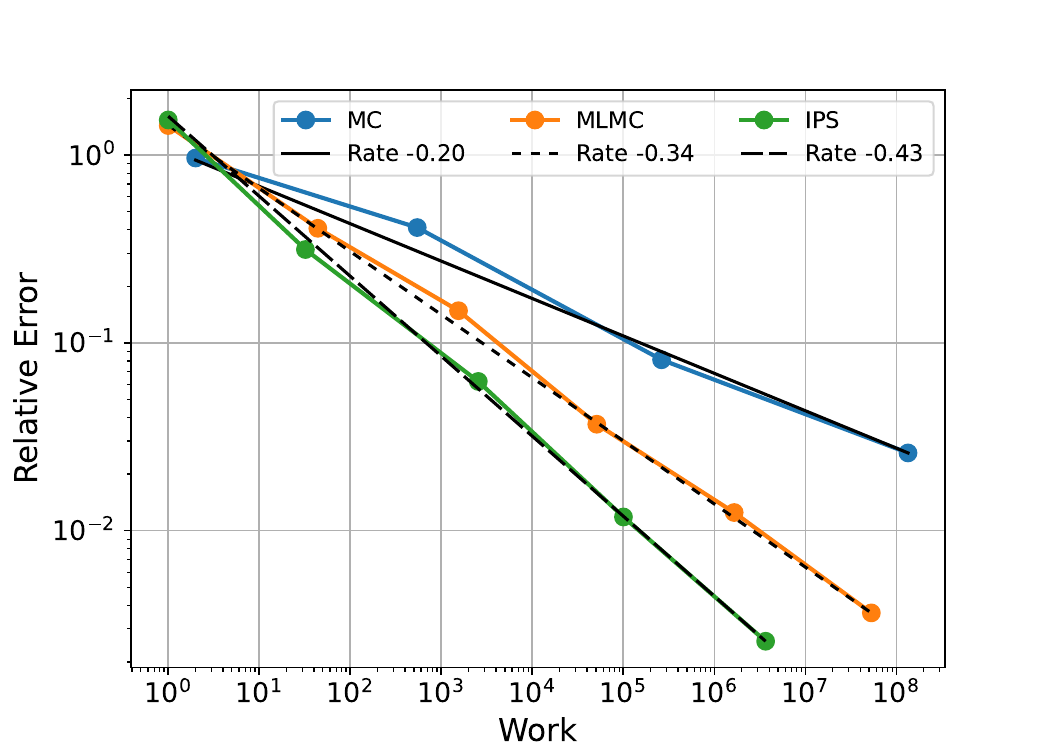}%
    \label{subfig:q2r5:b}%
  }\\
  \subfloat[$\vartheta=0.01$ and $\epsilon=0.05\vartheta$]{%
    \includegraphics[width=.4\linewidth]{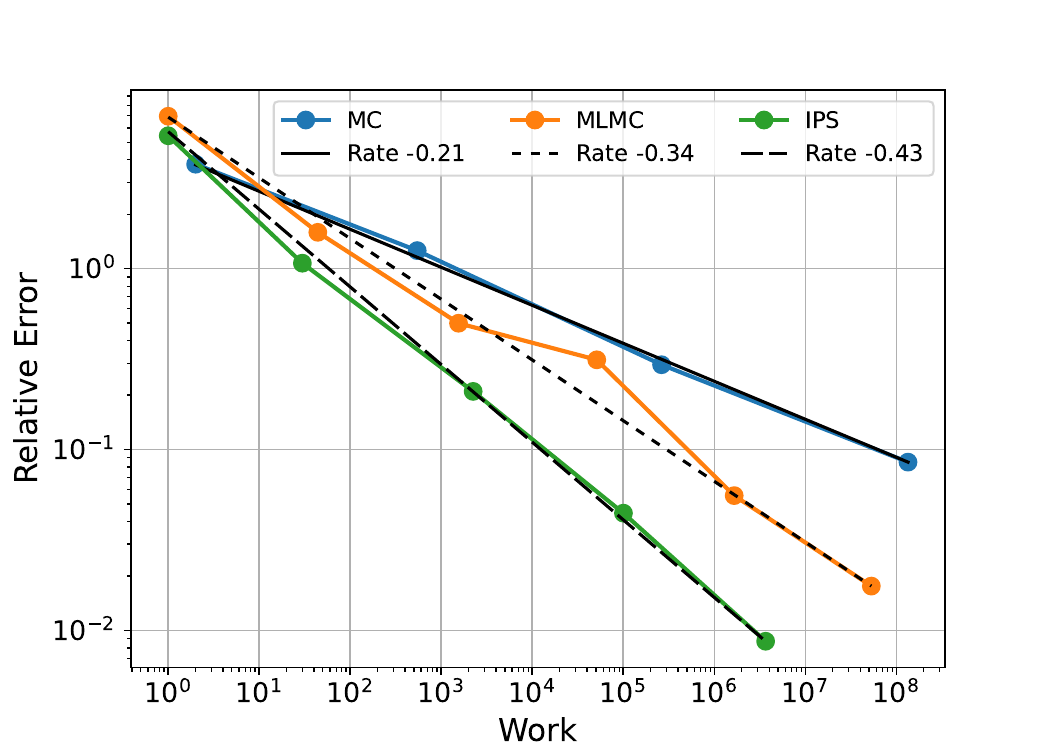}%
    \label{subfig:q2r5:c}%
  }\hfill
  \subfloat[$\vartheta=0.01$ and $\epsilon=0.01\vartheta$]{%
    \includegraphics[width=.4\linewidth]{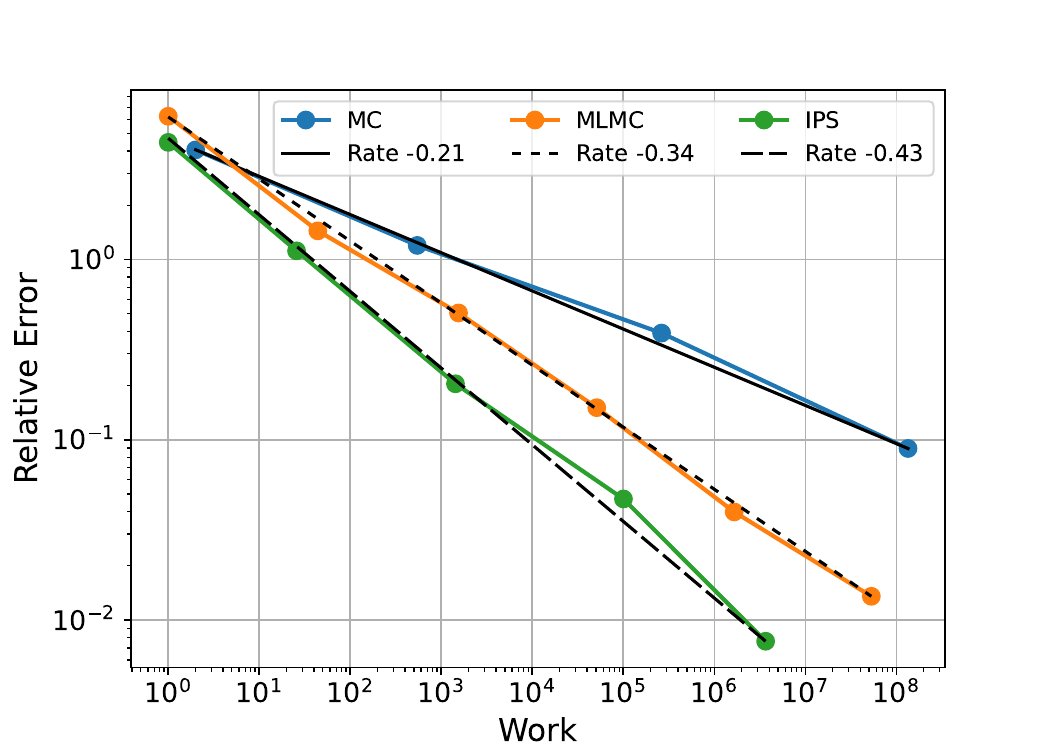}%
    \label{subfig:q2r5:d}%
  }
  \caption{Convergence for the Monte Carlo method, multilevel Monte Carlo method and multilevel interacting particle system method for $q=2$ and $r=5$ for different values of the cutoff $\vartheta$ and error constant $\epsilon$.}
  \label{fig:Simpleq2r5}
  % Run with:
  % python runExperiment.py --ips --mlmc --mc -L 8 -Y 0.1 -E 0.005 -r 4 -q 2 -g 0.5 --realizations 200 -f E.005_Y0.1_r4_q2_g0.5
  % python runExperiment.py --ips --mlmc --mc -L 8 -Y 0.1 -E 0.001 -r 4 -q 2 -g 0.5 --realizations 200 -f E.001_Y0.1_r4_q2_g0.5

  % python runExperiment.py --ips --mlmc --mc -L 8 -Y 0.01 -E 0.0005 -r 4 -q 2 -g 0.5 --realizations 200 -f E.0005_Y0.01_r4_q2_g0.5
  % python runExperiment.py --ips --mlmc --mc -L 8 -Y 0.01 -E 0.0001 -r 4 -q 2 -g 0.5 --realizations 200 -f E.0001_Y0.01_r4_q2_g0.5

\end{figure*}

Finally, we take $q=2$ and $r=5$ and repeat our experiment once again. The results for this last case are presented in \Cref{fig:Simpleq2r5} and for which we expect convergence rates of $\frac{q}{2q+r}=\frac{2}{9}$ for the MC method, $\frac{q}{q+r}=\frac{2}{7}$ for the MLMC method, and (similarly as for the previous experiment) a slightly worse convergence rate than $\frac{q}{r}=\frac{2}{5}$ for the MLIPS method (see the estimate in \Cref{thm:main} for this case). \Cref{fig:Simpleq2r5} displays the convergence results for this case, which follows the same trends as the experiments shown above, where the worsening of the convergence rates for the MLMC method and MLIPS method is stronger than before.

\subsubsection{Sampling close to the interface $\cG(y)=\vartheta$}
\begin{figure*}[ht!]
  \subfloat[$\ell=0$]{%
    \includegraphics[width=.35\linewidth]{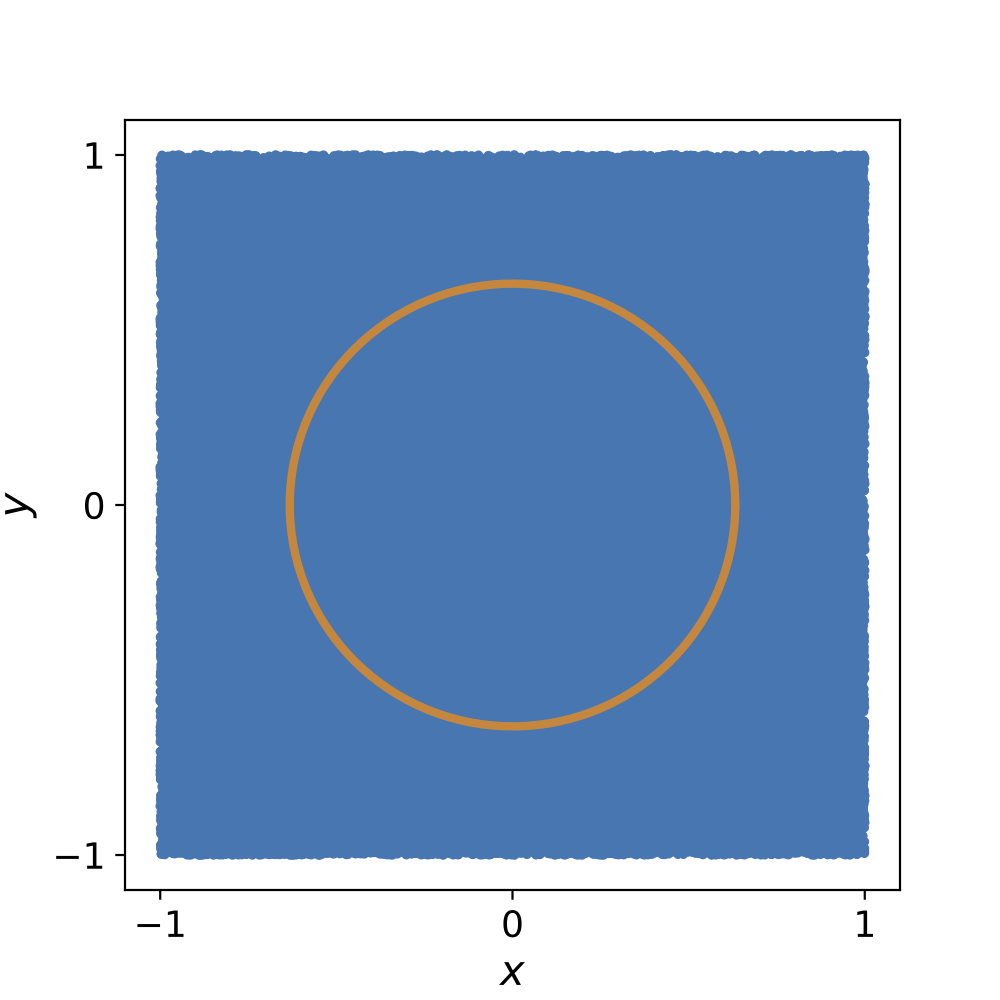}%
    \label{subfig:samples:a}%
  }\hfill
  \subfloat[$\ell=1$]{%
    \includegraphics[width=.35\linewidth]{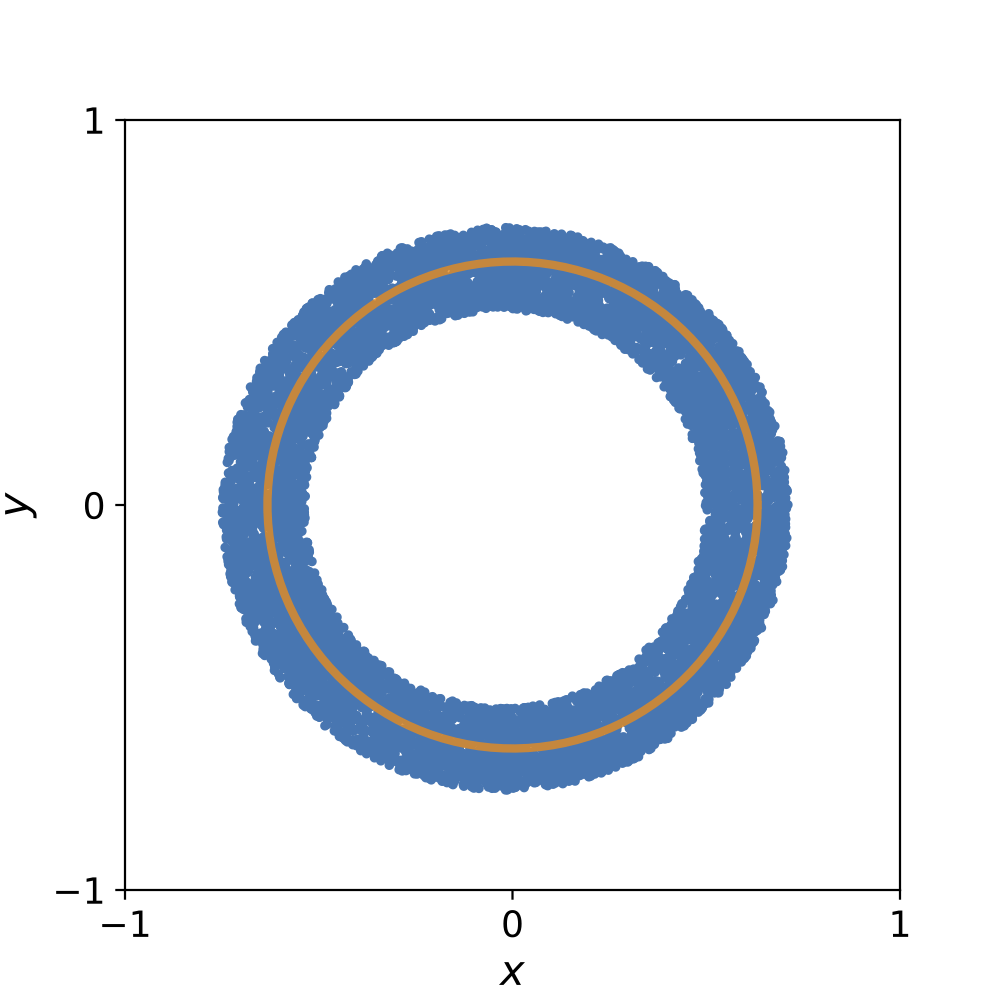}%
    \label{subfig:samples:b}%
  }\\
  \subfloat[$\ell=2$]{%
    \includegraphics[width=.35\linewidth]{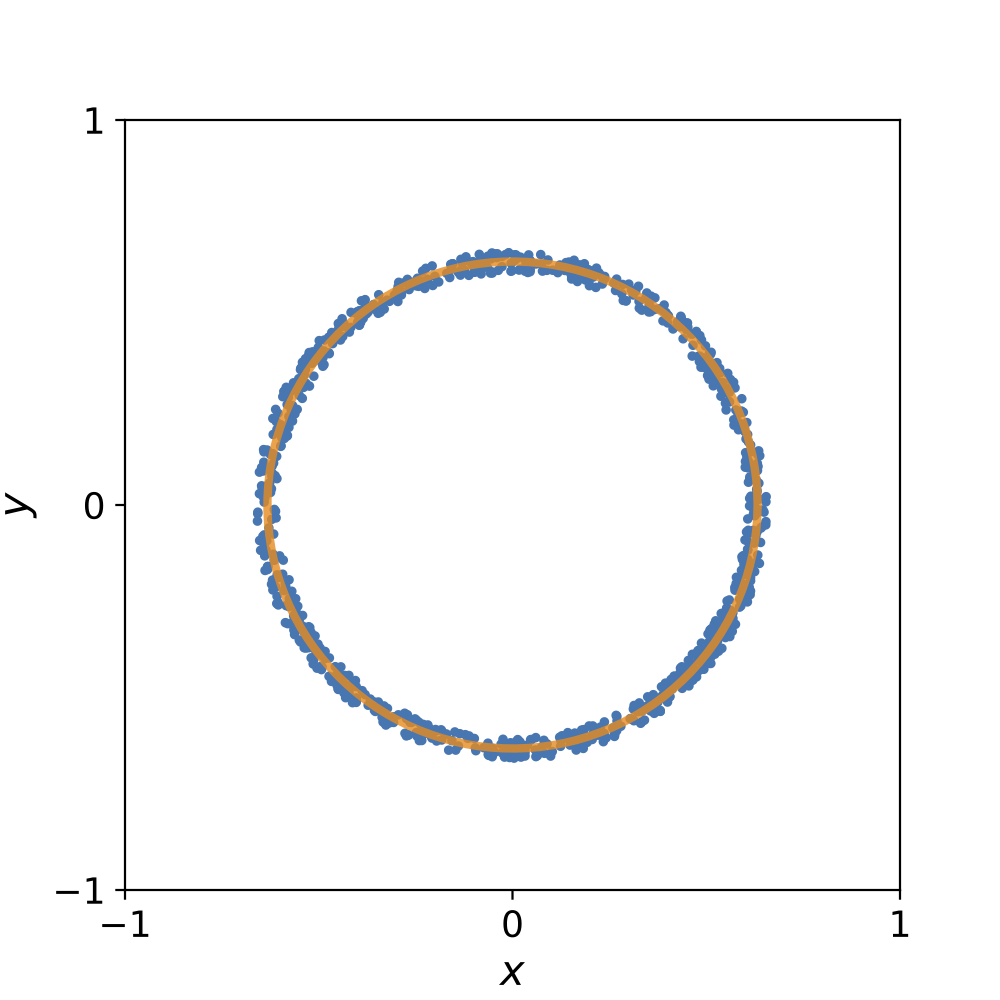}%
    \label{subfig:samples:c}%
  }\hfill
  \subfloat[$\ell=3$]{%
    \includegraphics[width=.35\linewidth]{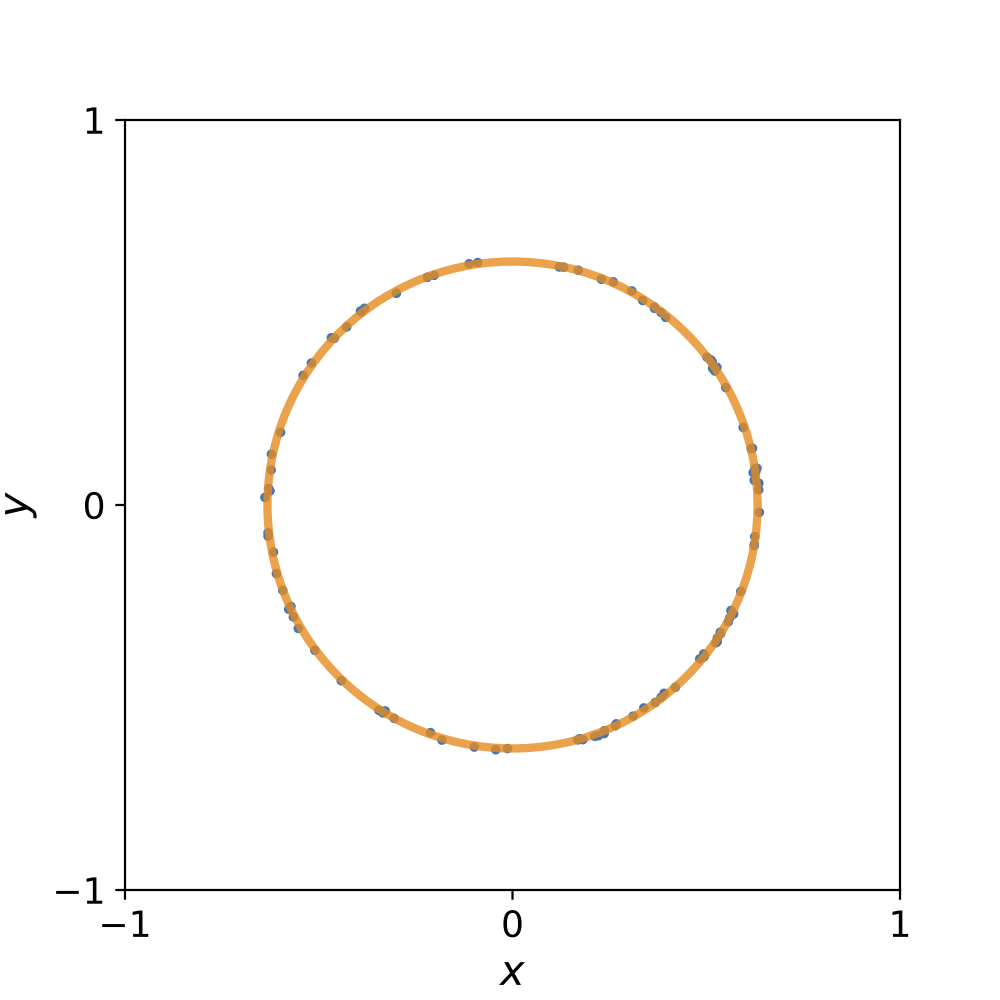}%
    \label{subfig:samples:d}%
  }
  \caption{Samples computed by the MLIPS method for different values $\ell\in\{0:3\}$. The centered circle (with radius $2\vartheta^{\frac{1}{2}}$)
    corresponds to the interface $\cG(y)=\vartheta$. We observe how the samples, represented by the points in the figures, accumulate around the interface
  with increasing precision.}
  \label{fig:Sample}
\end{figure*}

In this section we focus on presenting the second main property of the MLIPS method, namely, the generation of
samples close to the interface $\cG(y)=\vartheta$. \Cref{fig:Sample} displays, for the same parameters
used in \Cref{sssec:numqgr}--namely $q=2$ and $r=3$, $\epsilon=0.05\vartheta$ and $\vartheta=0.1$--the samples generated by the method for different levels for fixed $L=3$.
To further illustrate this property, we present a comparison against the method developed in \cite{elfversonMultilevelMonteCarlo2016}, where an adaptive scheme is proposed to enhance the efficiency of the Multilevel Monte Carlo by only computing high-resolution realizations of $\cG_\ell(y)$ when the sample points fall sufficiently close to the interface $\cG(y)=\vartheta$, which allows for increased sampling at (possibly) high resolutions. For further details we refer the interested reader to \cite{elfversonMultilevelMonteCarlo2016}. \Cref{tab:overallNumber} displays the average number of points that each algorithm recognizes as close to the boundary in each level $\ell\in\{1:3\}$ for 10 realizations and different values of $q$ and $r$, as well as the average computational cost required by each algorithm. For both tables we have taken $L=4$ and $\vartheta=0.1$, but \Cref{tab:eps5} considers $\epsilon=0.05\vartheta$ and \Cref{tab:eps1} considers $\epsilon=0.01\vartheta$. The number of samples for each algorithm was chosen so that the MLIPS generates $100$ samples at the last level and so that both methods required similar costs, simplifying the comparison of the number of generated samples. From the presented data, one can observe how the number of samples generated close to the limit interface shows very little variation for the MLIPS method, generating approximately $90$ points close to the interface for all cases (recall that the algorithm may generate repeated points during its call to \Cref{alg:sample}). The adaptive multilevel Monte Carlo method on the other hand, shows a strong dependence on the parameters $q$ and $\epsilon$, indicating that it finds it difficult to produce samples near the limit interface when the probability of the targeted area decreases, which occurs for lower values of $\epsilon$ and higher values of $q$ (which accelerate the decay of the area where higher resolution is required). Moreover, though the convergence properties presented in \Cref{thm:main} have some additional dependence on $L$ for the case $q\leq \frac{1}{2}r$ when compared to the convergence of the algorithm presented in \cite{elfversonMultilevelMonteCarlo2016}, this deterioration in the rate has little influence on our results, as seen in \Cref{fig:Simpleq2r4,fig:Simpleq2r5}.

\begin{table}[htbp]
  \centering
  \subfloat[$\vartheta=0.1$ and $\epsilon=0.05\vartheta$]{
    \begin{tabular}{l|l|llll|l}
      &   & \multicolumn{4}{c}{Samples per Level} & \\
      &   & $\ell=1$ & $\ell=2$   & $\ell=3$ & $\ell=4$ & Cost \\ \hline
      \multirow{2}{*}{$q=2$, $r=3$} & MLIPS & 100202.1 & 9462.4 & 933.0 & 90.2 & 4015100.0\\
      & MLAd & 24274.2 & 2161.7 & 177.4 & 17.2 & 4075659.2 \\ \hline
      \multirow{2}{*}{$q=1$, $r=3$} & MLIPS & 25201.0 & 3818.6 & 556.7 & 91.1 & 1633500.0 \\
      & MLAd & 7672.3 & 1362.5 & 226.9 & 42.8 & 1615924.0
    \end{tabular}
    \label{tab:eps5}
  }
  \hfill
  \subfloat[$\vartheta=0.1$ and $\epsilon=0.01\vartheta$]{
    \begin{tabular}{l|l|llll|l}
      &   & \multicolumn{4}{c}{Samples per Level} & \\
      &   & $\ell=1$ & $\ell=2$   & $\ell=3$ & $\ell=4$ & Cost \\ \hline
      \multirow{2}{*}{$q=2$, $r=3$} & MLIPS & 100174.2 & 9457.1 & 935.5 & 87.9 & 4015100.0\\
      & MLAd & 5544.0 & 494.3 & 36.6 & 3.4 & 4005643.0 \\ \hline
      \multirow{2}{*}{$q=1$, $r=3$} & MLIPS & 25208.2 & 3818.8 & 560.1 & 92.9 & 1633500.0 \\
      & MLAd & 2109.5 & 378.1 & 61.7 & 11.3 & 1625616.4
    \end{tabular}
    \label{tab:eps1}
  }
  \caption{Number of (unique) samples generated close to the interface $\cG(y)=\vartheta$ for the multilevel interacting particle system method (MLIPS) and the adaptive multilevel Monte Carlo (MLAd) method in \cite{elfversonMultilevelMonteCarlo2016}.}
  \label{tab:overallNumber}
\end{table}

\subsection{An applied example}
\begin{figure*}[ht!]
  \subfloat[$\ell=0$]{%
    \includegraphics[width=.4\linewidth]{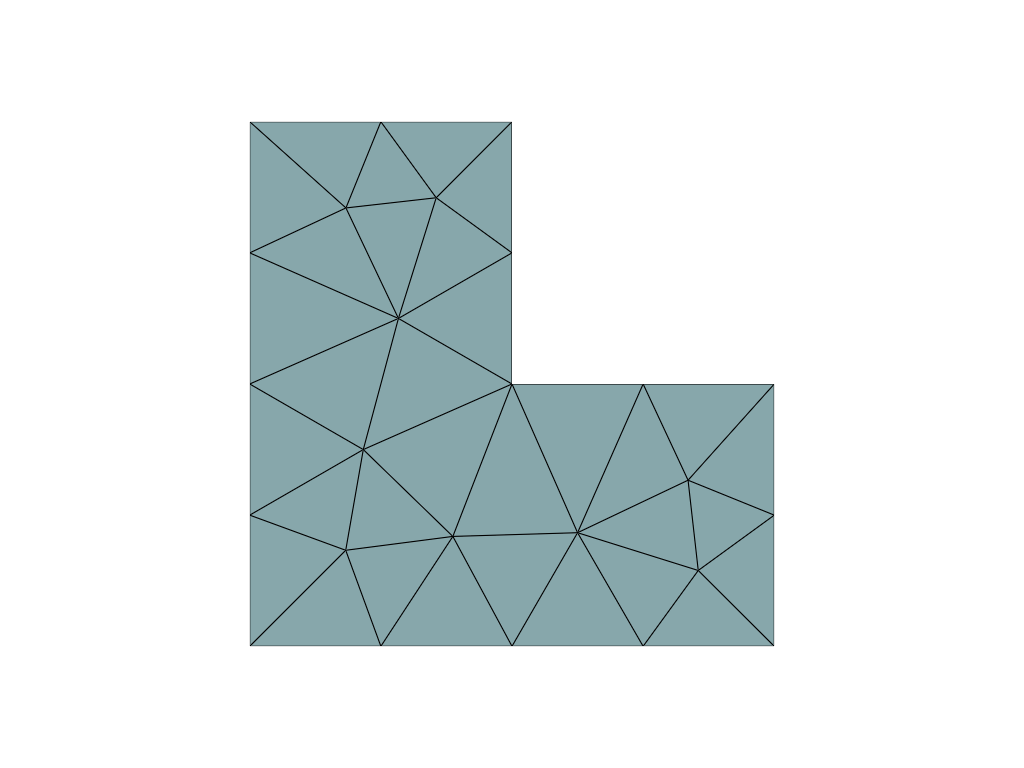}%
    \label{subfig:mesh:a}%
  }\hfill
  \subfloat[$\ell=1$]{%
    \includegraphics[width=.4\linewidth]{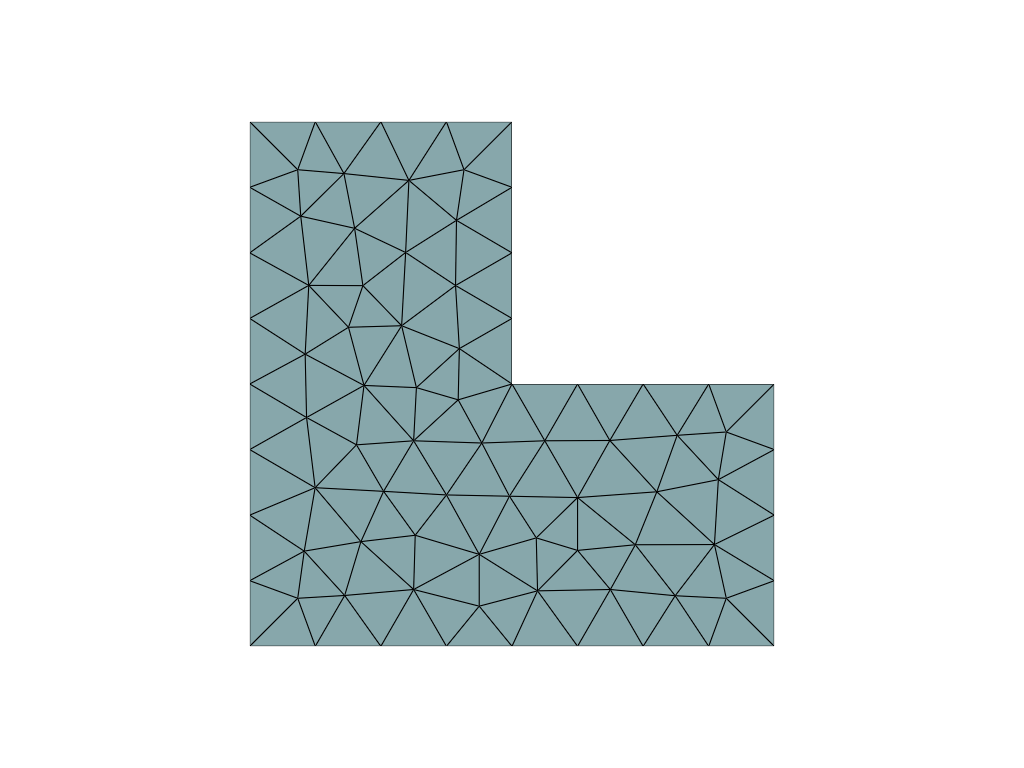}%
    \label{subfig:mesh:b}%
  }\\
  \subfloat[$\ell=2$]{%
    \includegraphics[width=.4\linewidth]{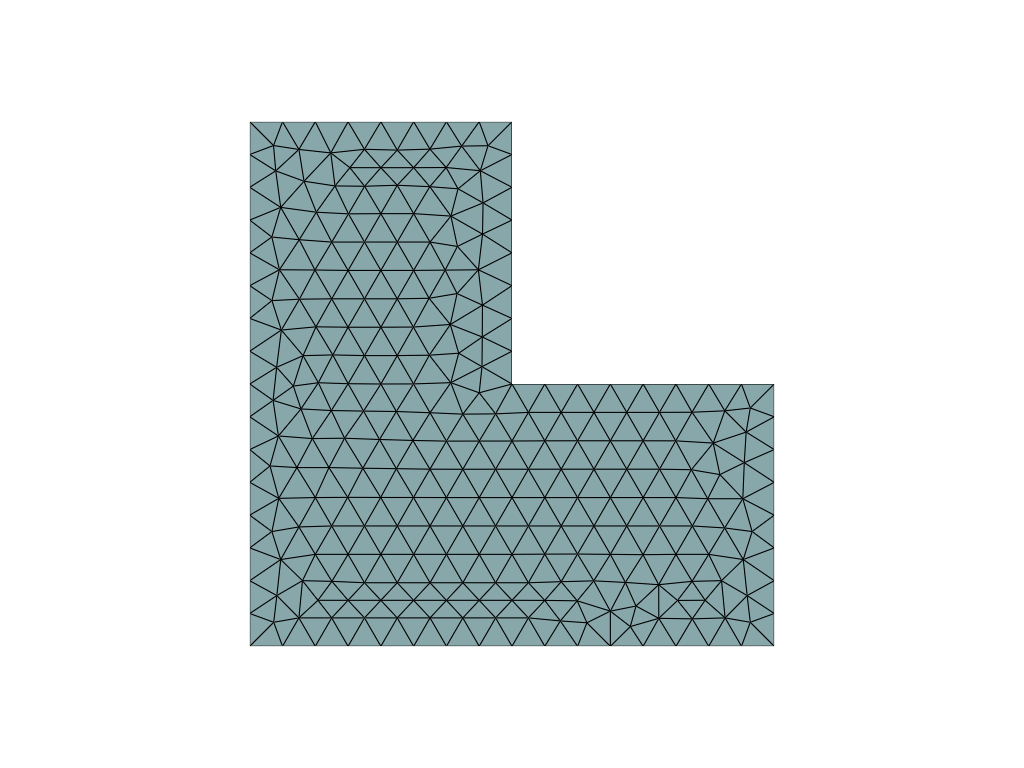}%
    \label{subfig:mesh:c}%
  }\hfill
  \subfloat[$\ell=3$]{%
    \includegraphics[width=.4\linewidth]{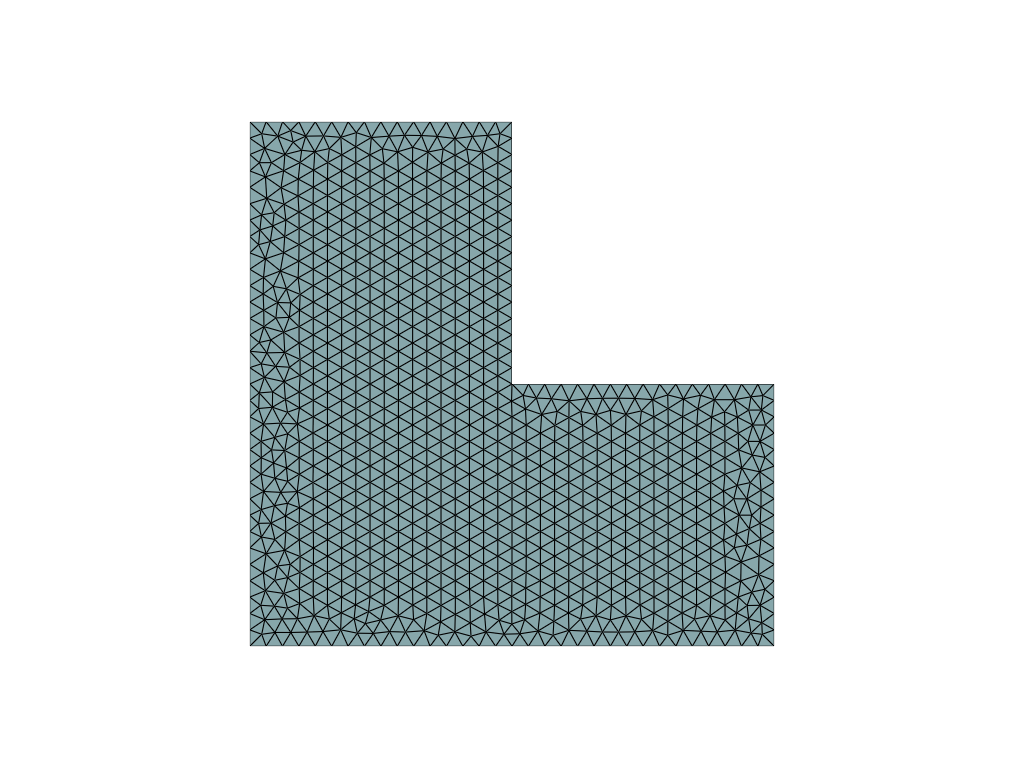}%
    \label{subfig:mesh:d}%
  }
  \caption{Meshes of $\Omega$ for $\ell\in\{0:3\}$.}
  \label{fig:meshs}
\end{figure*}

\begin{figure*}[ht!]
  \subfloat[]{%
    \includegraphics[width=.4\linewidth]{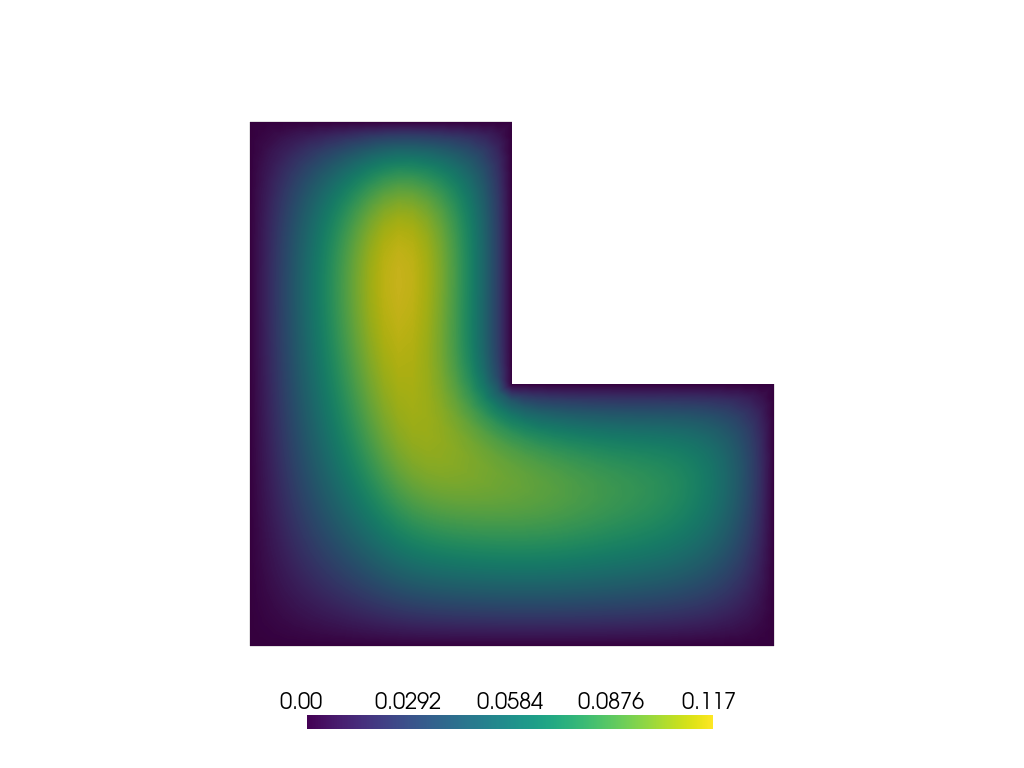}%
    \label{subfig:sol:a}%
  }\hfill
  \subfloat[]{%
    \includegraphics[width=.4\linewidth]{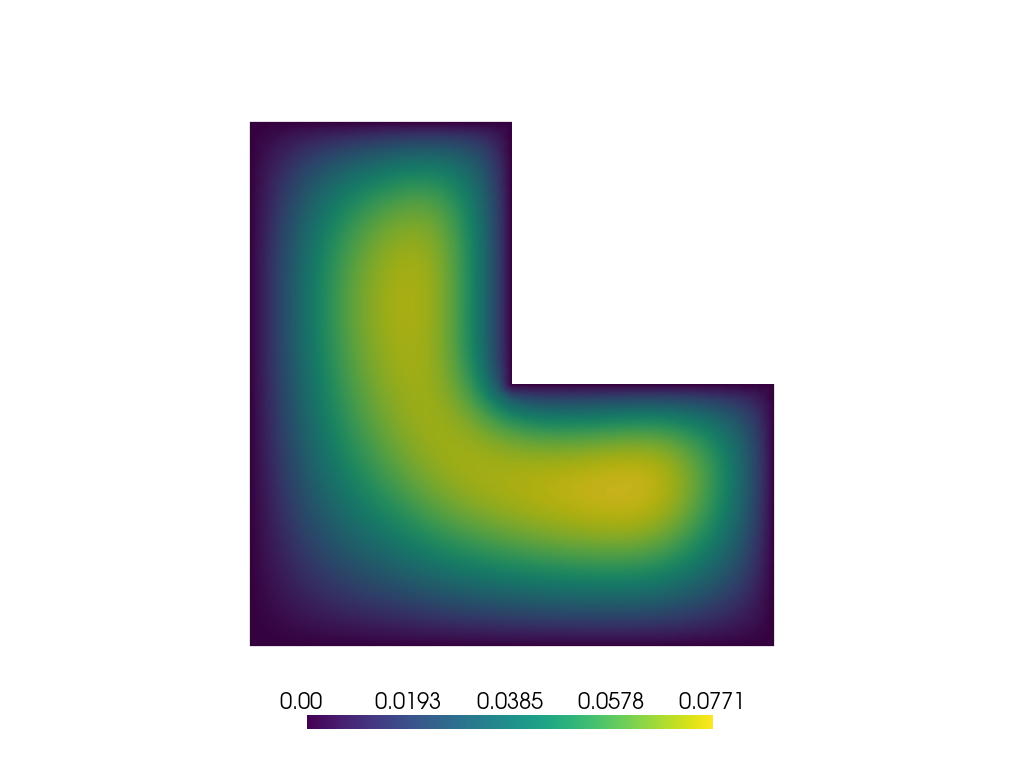}%
    \label{subfig:sol:b}%
  }\\
  \subfloat[]{%
    \includegraphics[width=.4\linewidth]{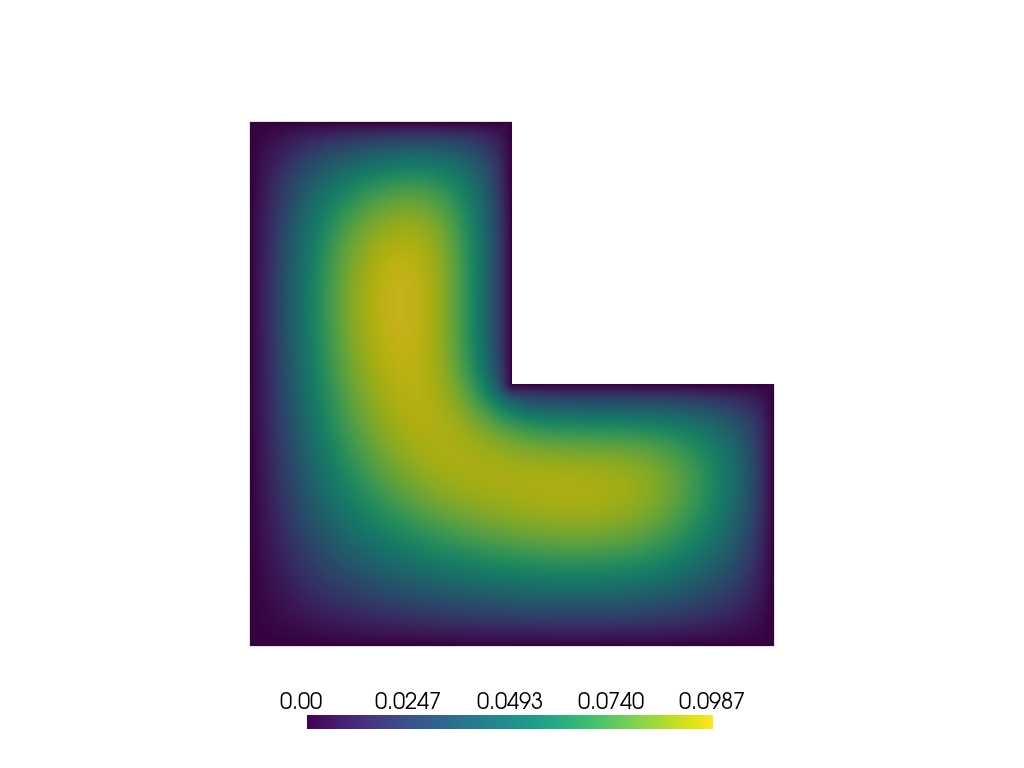}%
    \label{subfig:sol:c}%
  }\hfill
  \subfloat[]{%
    \includegraphics[width=.4\linewidth]{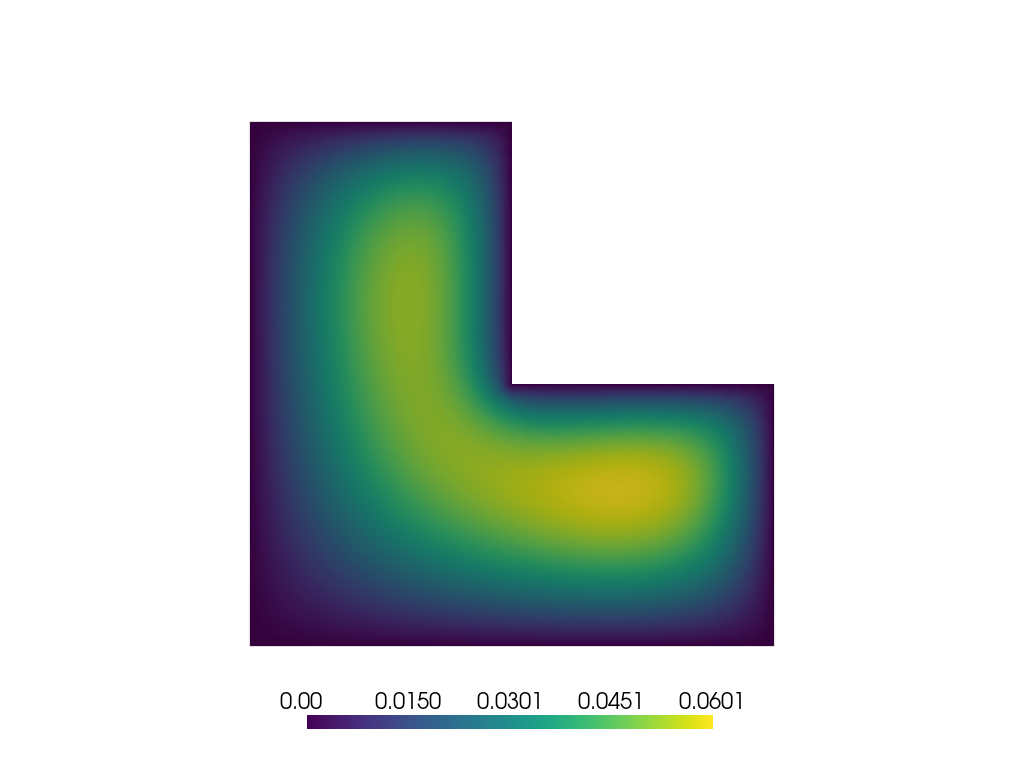}%
    \label{subfig:sol:d}%
  }
  \caption{Solutions to four different random parametric problems.}
  \label{fig:sols}
\end{figure*}

\begin{figure*}[ht!]
  \centering
  \subfloat[Convergence of the QoI]{%
    \includegraphics[width=.4\linewidth]{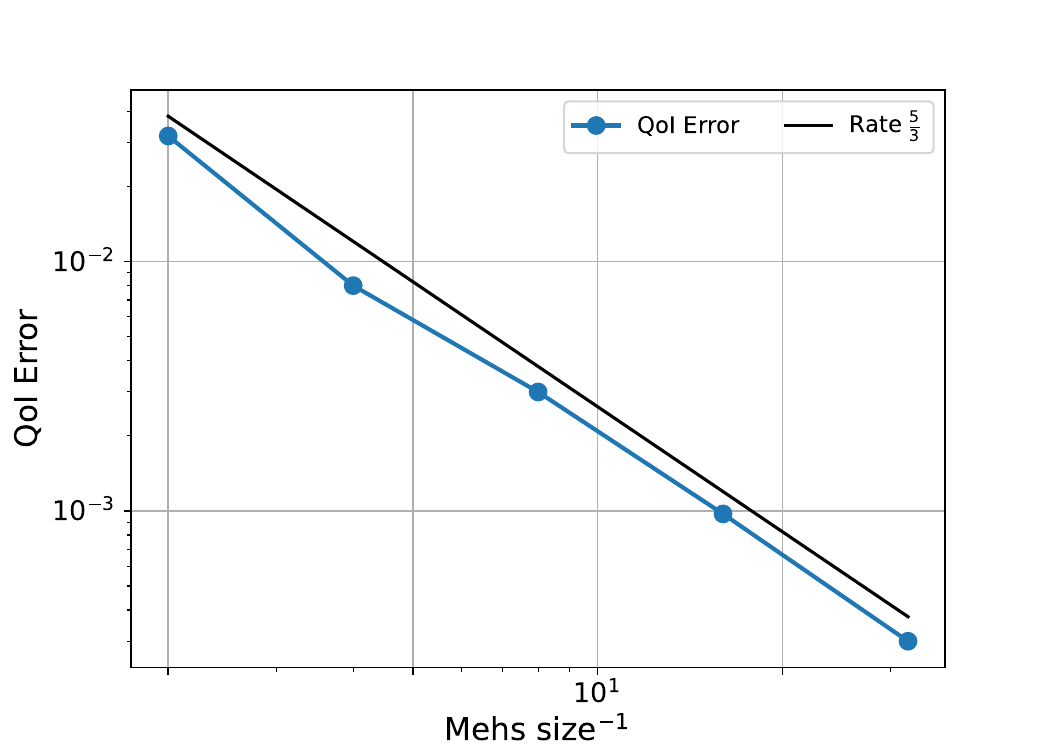}%
    \label{subfig:conv:a}%
  }\hfill
  \subfloat[Convergence of the failure probability]{%
    \includegraphics[width=.4\linewidth]{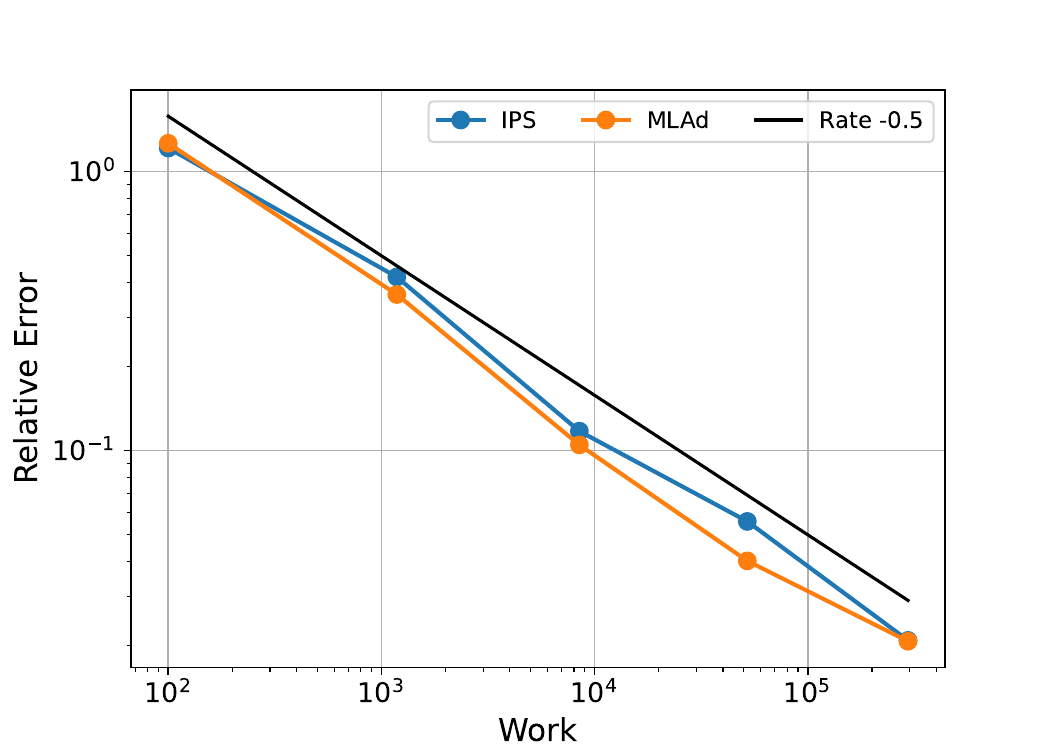}%
    \label{subfig:conv:b}%
  }
  \caption{Convergence of the QoI (for an arbitrary parameter, \Cref{subfig:conv:a}) and failure probability associated to the PDE in \Cref{eq:pde} (\Cref{subfig:conv:b}), for which
    $q=\tfrac{5}{3}$ and $r=2$. The reference value used to obtain \Cref{subfig:conv:a} was computed using an additionally refined mesh, while the results in \Cref{subfig:conv:b} was obtained by averaging $10$ realizations of the estimators for the given problem (which was enough to obtain the expected convergence rate) and compared to a reference value computed with the algorithm in \cite{elfversonMultilevelMonteCarlo2016} with an additional refinement level and an additional order of magnitude in the total work required for its computation.}
  \label{fig:convReal}
\end{figure*}

For our next example, and to once again display the properties of \Cref{alg:comp} in an applied setting, we consider a parametric PDE boundary value problem. For $m \in \IN$ we consider $\cP = [-1,1]^m$ equipped with the uniform measure. We further fix the L-shaped domain $\Omega:=(0,2)^2\setminus [1,2)^2$ and for any $y \in \cP$ we consider the problem of finding $v_y \in H^1(\Omega)$ such that
\begin{align}
  \label{eq:pde}
  \begin{aligned}
    -\nabla\cdot(d(y;x)\nabla v_y (x)) &= f(x)\quad\text{in}\quad\Omega,\\
    \gamma_d v_y(x) &= 0\quad\text{on}\quad\Gamma:=\partial\Omega,
  \end{aligned}
\end{align}
where $\gamma_d$ denotes the Dirichlet trace on the boundary $\Gamma$
and $d:[-1,1]^m\times\Omega\to\IR$ is a (strictly positive) parametric diffusion coefficient, which encodes
the randomness of the system in its first variable, and is given by
\begin{align*}
  d(y;x):=1+\frac{1}{4}\sum\limits_{i=1}^n\frac{y_i}{i}\left(1+\sin(\pi ix_1)\right).
\end{align*}
We will also fix
\begin{align*}
  f(x):=x_0x_1,
\end{align*} 
and consider the following QoI:
\begin{align*}
  \cG(y):=\int\limits_{\Omega} v_y(x)\,\dx.
\end{align*}
We approximate the solution of the PDE using the DOLFINx and UFL libraries in Python \cite{alnaesUnifiedFormLanguage2014,barattaDOLFINxNextGeneration2023} and Gmsh for mesh generation \cite{geuzaineGmsh3DFinite2009}. In particular, we consider a family of meshes such that the mesh at level $\ell$ has a meshsize given by $h_\ell=0.5^{\ell+1}$. We take $m=8$ random parameters and, given the reentrant corner in the geometry of $\Omega$, expect a convergence rate of $\tfrac{5}{3}$ for the chosen quantity of interest due to the regularity of the solution to the PDE in \cref{eq:pde} (we refer to \cite{caiFiniteElementMethod2001,costabelANALYTICREGULARITYLINEAR2012} and references therein for more details regarding the regularity of solutions in non-convex polygonal domains). \Cref{subfig:conv:a} displays the convergence of the QoI to a reference value computed with an additionally refined mesh, while \Cref{fig:meshs,fig:sols} show, respectively, the first four meshes of $\Omega$ ($\ell\in\{0:3\}$) and four solutions of the PDE in \cref{eq:pde} for four arbitrary parameters. Given the properties of the problem, we take $r=2$ (to display a computational work growing proportionally with the degrees of freedom) and $q=\tfrac{5}{3}$. The constant $C_\cG$ is estimated numerically. \Cref{tab:meshs} shows the mesh information for different values of $\ell$, where we see how the degrees of freedom roughly increase by a factor of $4$ as the meshsize is halved. Finally, \Cref{subfig:conv:b} shows the attained convergence rates for the computation of the probability through the proposed algorithm and the adaptive algorithm introduced in \cite{elfversonMultilevelMonteCarlo2016}, showing how both algorithms achieve the predicted rate of $0.5$ with respect to the total computational effort in the range $q\in(\tfrac{r}{2},r)$ (in accordance with \Cref{thm:main}). The reference value used to obtain the results in \Cref{subfig:conv:b} was computed through the adaptive multilevel Monte Carlo algorithm in \cite{elfversonMultilevelMonteCarlo2016}.

\begin{table}[]
  \centering
  \begin{tabular}{l|llllll}
    &\multicolumn{6}{c}{Mesh Information  per Level} \\
    Level &$\ell=0$ & $\ell=1$  & $\ell=2$ & $\ell=3$ & $\ell=4$  & $\ell=5$ \\ \hline
    Meshsize &  0.5 & 0.25 & 0.125 & 0.0625 & 0.03125 & 0.015625\\
    DoFs &  32 & 128 & 382 & 1816 & 7170 & 28490 \\\hline
  \end{tabular}
  \caption{Number of Mesh Elements per level.}
  \label{tab:meshs}
\end{table}

\section{Conclusions}
\label{sec:conclusions}
We have proposed a novel multilevel method for the computation of failure probabilities in the context of computational uncertainty quantification and, more specifically, with Partial Differential Equations depending on random parameters in mind. The algorithm leverages a sequential sampling scheme in order to sample the multilevel contributions only where meaningful, i.e., where there is a positive probability of the contribution being different from $0$. This allows the proposed algorithm to achieve more efficient convergence rates than those achieved by classical Monte Carlo and multilevel Monte Carlo methods \cite{elfversonMultilevelMonteCarlo2016,haji-aliAdaptiveMultilevelMonte2022}, while achieving comparable convergence rates to the current best alternative introduced in \cite{elfversonMultilevelMonteCarlo2016}. Moreover, the algorithm naturally generates samples close to the relevant limit-state surface, allowing later usage for possible applications and extensions through reduced order models \cite{dasguptaFailureProbabilityEstimation2019} and computation of sensitivities of the failure probability \cite{papaioannouReliabilitySensitivityEstimation2018}. To the best of our knowledge, this is the only algorithm which has leveraged such a sequential sampling scheme (based on the theory in \cite{cerouSequentialMonteCarlo2012,cerouNonasymptoticTheoremUnnormalized2011,DelMoralFeynmanKac2004,delmoralGenealogicalParticleAnalysis2005,moralMultilevelSequentialMonte2017} and references therein) to the multilevel computation of failure probabilities in such a way. Numerical examples then display the convergence properties of the algorithm both for a simple fabricated example (allowing for easy extensive experimentation) and an applied example concerning the Laplace equation on an L-shaped domain. Future research directions include the usage of the proposed algorithm in the context of robust uncertainty quantification, enhancing the algorithm with reduced order models (in comparison to, for example, the work in \cite{aylwinUncertaintyQuantificationDiffraction2025}, where a reduced basis is built on the whole sample space to approximate a failure probability) and possible generalizations in the direction of rare event estimation.

\bibliographystyle{plain}
\bibliography{references,manualReferences}
\appendix
\section{Number of Samples per Level for the MLIPS Approximation} 
\label{apdx:samp}
We briefly describe how to determine the number of samples per level $\cN_L$ so as to obtain
the results displayed in \cref{thm:main}. As mentioned before, our proof is based on that of
\cite[Thm.~1]{cliffeMultilevelMonteCarlo2011}. We fix $L\in\IN$ and recall that we have the
following bound for the error of the MLIPS approximation:
\begin{align*}
  \IE\left(\vert\IE(\IH_{\cG})-\mlmk{\cN_L}{\IH_{\cG}}\vert^2\right)^{\frac{1}{2}}
  &\lesssim\alpha^{qL}+
    \sum\limits_{\ell=0}^L(\ell^{\frac{1}{2}}+1)\alpha^{q\ell}N_{\ell}^{-\frac{1}{2}},
\end{align*}
while for the required computational work it holds that
\begin{align*}
  \work{\mlmk{\cN_L}{\IH_\cG}}\lesssim\sum\limits_{\ell=0}^{L}N_\ell(\cW_\ell+\cW_{\ell-1})\lesssim\sum\limits_{\ell=0}^{L}N_\ell\alpha^{-r\ell}.
\end{align*}
We distinguish three different cases, corresponding to $q>\frac{1}{2}r$, $q=\frac{1}{2}r$, and
$q<\frac{1}{2}r$.
\subsection{Case $q>\frac{1}{2}r$}~\\
We choose $N_\ell=\Theta((\ell^{\frac{1}{2}}+1)^2\alpha^{-2qL}\alpha^{\frac{2}{3}(q+r)\ell})$. Then, it follows for the error of the method that
\begin{align*}
  \IE\left(\vert\IE(\IH_{\cG})-\mlmk{\cN_L}{\IH_{\cG}}\vert^2\right)^{\frac{1}{2}}
  &\lesssim\alpha^{qL}+
    \sum\limits_{\ell=0}^L(\ell^{\frac{1}{2}}+1)\alpha^{q\ell}N_{\ell}^{-\frac{1}{2}},\\
  &\lesssim\alpha^{qL}+
    \alpha^{qL}\sum\limits_{\ell=0}^L\alpha^{q\ell}\alpha^{-\frac{1}{3}(q+r)\ell}\\
  &=\alpha^{qL}+
    \alpha^{qL}\sum\limits_{\ell=0}^L\alpha^{\frac{1}{3}(2q-r)\ell}
    \lesssim\alpha^{qL},
\end{align*}
where the last equation follows from $q>\frac{1}{2}r$, which implies that the last sum is bounded independently of $L\in\IN$.
For the computational work, it holds that
\begin{align*}
  \work{\mlmk{\cN_L}{\IH_\cG}}\lesssim\sum\limits_{\ell=0}^{L}N_\ell(\cW_\ell+\cW_{\ell-1})
  &\lesssim\sum\limits_{\ell=0}^{L}N_\ell\alpha^{-r\ell}\\
  &\lesssim\alpha^{-2qL}\sum\limits_{\ell=0}^{L}(\ell^{\frac{1}{2}}+1)^2\alpha^{-r\ell}\alpha^{\frac{2}{3}(q+r)\ell}\\
  &=\alpha^{-2qL}\sum\limits_{\ell=0}^{L}(\ell^{\frac{1}{2}}+1)^2\alpha^{\frac{1}{3}(2q-r)\ell}
    \lesssim\alpha^{-2qL},
\end{align*}
where the last equation follows, just as before, from the fact that the last sum is bounded independently of $L\in\IN$.
\subsection{Case $q<\frac{1}{2}r$}~\\
We choose $N_\ell=\Theta((\ell^{\frac{1}{2}}+1)^2\alpha^{-2qL}\alpha^{\frac{2}{3}(2q-r)L}\alpha^{\frac{2
  }{3}(q+r)\ell})$ just as for the case $q>\frac{1}{2}r$, which gives
\begin{align*}
  \IE\left(\vert\IE(\IH_{\cG})-\mlmk{\cN_L}{\IH_{\cG}}\vert^2\right)^{\frac{1}{2}}
  &\lesssim\alpha^{qL}+
    \alpha^{qL}\alpha^{\frac{1}{3}(r-2q)L}\sum\limits_{\ell=0}^L\alpha^{\frac{1}{3}(2q-r)\ell}\\
  &=\alpha^{qL}+
    \alpha^{qL}\sum\limits_{\ell=0}^L\alpha^{\frac{1}{3}(2q-r)(\ell-L)}\\
  &=\alpha^{qL}+
    \alpha^{qL}\sum\limits_{\ell=0}^L\alpha^{\frac{1}{3}(r-2q)\ell}
    \lesssim \alpha^{qL}
\end{align*}
where again the last sum is bounded independently of $L\in\IN$. For the computational work, proceeding as before yields
\begin{align*}
  \work{\mlmk{\cN_L}{\IH_\cG}}
  &\lesssim \alpha^{-2qL}\alpha^{\frac{2}{3}(2q-r)L}\sum\limits_{\ell=0}^{L}(\ell^{\frac{1}{2}}+1)^2\alpha^{\frac{1}{3}(2q-r)\ell}\\
  &\leq \alpha^{-2qL}\alpha^{\frac{2}{3}(2q-r)L}(L^{\frac{1}{2}}+1)^2\sum\limits_{\ell=0}^{L}\alpha^{\frac{1}{3}(2q-r)\ell}\\
  &= \alpha^{-2qL}\alpha^{(2q-r)L}(L^{\frac{1}{2}}+1)^2\sum\limits_{\ell=0}^{L}\alpha^{\frac{1}{3}(r-2q)\ell}\lesssim L\alpha^{-rL}.
\end{align*}
\subsection{Case $q=\frac{1}{2}r$}~\\
We choose $N_\ell=\Theta(L^2(\ell^{\frac{1}{2}}+1)^2\alpha^{-2qL}\alpha^{2q\ell})$, from where it follows that
\begin{align*}
  \IE\left(\vert\IE(\IH_{\cG})-\mlmk{\cN_L}{\IH_{\cG}}\vert^2\right)^{\frac{1}{2}}
  &\lesssim\alpha^{qL}+
    \alpha^{qL}L^{-1}\sum\limits_{\ell=0}^L\alpha^{\frac{1}{3}(2q-r)\ell}\\
  &=\alpha^{qL}+
    \alpha^{qL}L^{-1}\sum\limits_{\ell=0}^L1\lesssim \alpha^{qL}.
\end{align*}
Finally, we have that the computational work may be bounded as
\begin{align*}
  \work{\mlmk{\cN_L}{\IH_\cG}}
  &\lesssim \alpha^{-2qL}L^2\sum\limits_{\ell=0}^{L}(\ell^{\frac{1}{2}}+1)^2\alpha^{(2q-r)\ell}\\
  &\lesssim \alpha^{-2qL}L^2\sum\limits_{\ell=0}^{L}(\ell^{\frac{1}{2}}+1)^2\lesssim L^4
    \alpha^{-rL}.
\end{align*}

\end{document}